\documentclass[a4paper,12pt]{amsart}
\usepackage[utf8]{inputenc}
\usepackage[english]{babel} 

\usepackage[T1]{fontenc}

\usepackage{amsthm}
\usepackage{amsmath}
\usepackage{amsfonts}
\usepackage{mathtools}
\usepackage[colorlinks]{hyperref}
\usepackage{enumerate}

\def\bdi{\begin{diagram}}
\def\edi{\end{diagram}}

\newtheorem{thm}{Theorem}[section]
\newtheorem{cor}[thm]{Corollary}
\newtheorem{lem}[thm]{Lemma}
\newtheorem{prop}[thm]{Proposition}
\theoremstyle{definition}
\newtheorem{defi}[thm]{Definition}
\newtheorem{defis}[thm]{Definitions}
\newtheorem{conj}[thm]{Conjecture}
\newtheorem{conv}[thm]{Convention}
\newtheorem{nota}[thm]{Notation}
\newtheorem{rem}[thm]{Remark}
\newtheorem{rems}[thm]{Remarks}
\newtheorem{exa}[thm]{Example}
\newtheorem{exas}[thm]{Examples}
\newtheorem{prob}[thm]{Problem}
\newtheorem{probs}[thm]{Problems}
\newtheorem{ques}[thm]{Question}
\newtheorem{sett}[thm]{Setting}
\newtheorem{sit}[thm]{}

\newcommand{\Spec}{\operatorname{{\rm Spec}}}

\newcommand{\Hol}{ \operatorname{{\rm Hol}}}
\newcommand{\Ker}{ \operatorname{{\rm Ker}}}

\newcommand{\Aut}{ \operatorname{{\rm Aut}}}
\newcommand{\Mat}{\operatorname{{\rm Mat}}}
\newcommand{\LND}{\operatorname{{\rm LND}}}
\newcommand{\GL}{\operatorname{{\rm GL}}}
\newcommand{\SL}{ \operatorname{{\rm SL}}}

\def\pr{\mathop{\rm pr}}

\def\lto{\longrightarrow}

\renewcommand{\epsilon}{\varepsilon}

\def\and{\quad\mbox{and}\quad}

\renewcommand{\div}{ \operatorname{\rm div}}

\newcommand{\C}{\ensuremath{\mathbb{C}}}

\newcommand{\Q}{\ensuremath{\mathbb{Q}}}

\newcommand{\N}{\ensuremath{\mathbb{N}}}
\newcommand{\G}{\ensuremath{\mathbb{G}}}

\newcommand{\bK}{{\ensuremath{\rm \bf K}}}
\newcommand{\bk}{{\ensuremath{\rm \bf k}}}

\newcommand{\hI}{{\hat I}}
\newcommand{\hA}{{\hat A}}

\newcommand{\tB}{{\tilde B}}
\newcommand{\tC}{{\tilde C}}

\newcommand{\tH}{{\tilde H}}

\newcommand{\tY}{{\tilde Y}}

\newcommand{\bC}{{\bar C}}

\newcommand{\bZ}{{\bar Z}}

\newcommand{\cW}{{\ensuremath{\mathcal{W}}}}

\def\fl{{\mathfrak l}}

\newcommand{\cB}{{\ensuremath{\mathcal{B}}}}

\newcommand{\cF}{{\ensuremath{\mathcal{F}}}}
\newcommand{\cG}{{\ensuremath{\mathcal{G}}}}

\newcommand{\cE}{{\ensuremath{\mathcal{E}}}}
\newcommand{\cA}{{\ensuremath{\mathcal{A}}}}
\newcommand{\cJ}{{\ensuremath{\mathcal{J}}}}
\newcommand{\cO}{{\ensuremath{\mathcal{O}}}}
\newcommand{\cC}{{\ensuremath{\mathcal{C}}}}
\newcommand{\cD}{{\ensuremath{\mathcal{D}}}}

\newcommand{\cP}{{\ensuremath{\mathcal{P}}}}

\newcommand{\cR}{{\ensuremath{\mathcal{R}}}}
\newcommand{\cN}{{\ensuremath{\mathcal{N}}}}

\newcommand{\cX}{{\ensuremath{\mathcal{X}}}}

\newcommand{\cV}{{\ensuremath{\mathcal{V}}}}

\newcommand{\cZ}{{\ensuremath{\mathcal{Z}}}}
\newcommand{\p}{\partial}

\newcommand{\id}{{\rm id}}

\renewcommand{\rho}{\varrho}

\def\bals#1\eals{\begin{align*}#1\end{align*}}
\def\bal#1\eal{\begin{align}#1\end{align}}

\def\SAut{\mathop{\rm SAut}}
\def\SL{\mathop{\rm SL}}
\def\EE{\mathop{\rm E}}
\def\FF{\mathop{\rm F}}

\def\A{{\mathbb A}}

\def\CC{{\mathbb C}}

\def\PP{{\mathbb P}}

\def\fq{{\mathfrak q}}
\def\fp{{\mathfrak p}}

\def\fs{{\mathfrak s}}

\def\dd{ {\rm d }}

\renewcommand{\phi}{\varphi}

\newcommand{\bnum}{\begin{enumerate}}
\newcommand{\enum}{\end{enumerate}}

\newcommand{\brem}{\begin{rem}}
\newcommand{\brems}{\begin{rems}}
\newcommand{\erem}{\end{rem}}
\newcommand{\erems}{\end{rems}}
\newcommand{\bprob}{\begin{prob}}
\newcommand{\eprob}{\end{prob}}
\newcommand{\bprobs}{\begin{probs}}
\newcommand{\eprobs}{\end{probs}}
\newcommand{\bques}{\begin{ques}}
\newcommand{\eques}{\end{ques}}
\newcommand{\bexa}{\begin{exa}}
\newcommand{\bexas}{\begin{exas}}
\newcommand{\eexa}{\end{exa}}
\newcommand{\eexas}{\end{exas}}
\newcommand{\bdefi}{\begin{defi}}
\newcommand{\edefi}{\end{defi}}
\newcommand{\bdefis}{\begin{defis}}
\newcommand{\edefis}{\end{defis}}
\newcommand{\bcor}{\begin{cor}}
\newcommand{\ecor}{\end{cor}}
\newcommand{\blem}{\begin{lem}}
\newcommand{\elem}{\end{lem}}
\newcommand{\bconv}{\begin{conv}}
\newcommand{\econv}{\end{conv}}
\newcommand{\bconj}{\begin{conj}}
\newcommand{\econj}{\end{conj}}
\newcommand{\bprop}{\begin{prop}}
\newcommand{\eprop}{\end{prop}}
\newcommand{\bthm}{\begin{thm}}
\newcommand{\ethm}{\end{thm}}
\newcommand{\bnota}{\begin{nota}}
\newcommand{\enota}{\end{nota}}
\newcommand{\bsit}{\begin{sit}}
\newcommand{\esit}{\end{sit}}
\newcommand{\be}{\begin{equation}}
\newcommand{\ee}{\end{equation}}
\newcommand{\bproof}{\begin{proof}}
\newcommand{\eproof}{\end{proof}}
\newcommand{\bsett}{\begin{sett}}
\newcommand{\esett}{\end{sett}}
\def\ba{\begin{array}}
\def\ea{\end{array}}

\begin{document}
\title{On automorphisms of flexible varieties}

\author{
S.\ Kaliman and
D.\ Udumyan}
\address{Department of Mathematics,
University of Miami, Coral Gables, FL 33124, USA}
\email{kaliman@math.miami.edu}
\address{Department of Mathematics,
University of Miami, Coral Gables, FL 33124, USA}
\email{mathudum@gmail.com}


\date{\today}
\maketitle

\begin{abstract} Let $X$ be a flexible variety and $\varphi : C_1 \to C_2$ be an isomorphism
of closed one-dimensional subschemes of $X$. We develop several criteria which guarantee that $\varphi$ extends
to an automorphism of $X$. \end{abstract}


{\footnotesize \tableofcontents}

\section*{Introduction} 

Every algebraic variety $X$ in this paper is considered over an algebraically closed field $\bk$ of characteristic 0. 
We study the following question: 

\noindent {\em when can an isomorphism
 $f :Y_1 \to Y_2$ of subschemes of $X$  be extended to an automorphism of $X$?}
 
 For $X\simeq \A^n$ and $Y_i\simeq \A^1$  the answer to this question is known to be always  positive in all cases
 but $n=3$ (see \cite{AbMo}, \cite{Su},  \cite{Cr}, \cite{Je}). For $n=3$ the answer is unknown but if $\bk =\C$, then
 $f$ can be extended to a holomorphic automorphism of $\C^3$ \cite{Ka92}. Furthermore, very general sufficient conditions 
 for extension of isomorphisms of closed subvarieties of $\A^n$ to automorphisms of $\A^n$ where found in \cite{Ka91} and \cite{Sr}
 (in a weaker form they were already present in \cite{Je}).
 For $X\ne \A^n$ the first break-through is due
 to van Santen (formerly Stampfli) who proved that every isomorphism of smooth polynomial curves in $X\simeq \SL_n(\C), \, n\geq 3$ 
 is extendable to an automorphism of $X$ \cite{St}. Later jointly with Feller he showed that the same is true when $Y_i \simeq \C$
 and $X$ is any complex linear algebraic group  different from $\C^3$, $\SL_2(\C)$ or ${\rm PSL}_2 (\C)$ \cite{FS}.
 Some extra cases with a positive answer  to the extension problem were found in \cite{Ka20}. Say, the answer is positive when
 $X$ is a nonzero fiber  of a non-degenerate quadric form in $\A^m$ (we call such $X$ an $(m-1)$-sphere over $\bk$)
and $Y_i$ is a smooth closed subvariety with $\dim Y_i \leq \frac{m}{3}-1$.

Starting from dimension 2 affine spaces, linear algebraic groups without nontrivial characters and  spheres over $\bk$ are examples
of so-called flexible varieties  (recall that a normal quasi-affine variety $X$ is flexible if $\SAut (X)$ acts transitively on the smooth part
$X_{\rm reg}$ of $X$ where $\SAut (X)$ is
the subgroup of the group $\Aut (X)$ of algebraic automorphisms of $X$ generated by one-parameter unipotent subgroups). 
For a smooth flexible variety $X$  every isomorphism of zero-dimensional subvarieties is extendable to an automorphism \cite{AFKKZ}. Furthermore,
recall that for a closed subvariety  $Y$ of $X$ with  defining ideal $I$ in the ring $A=\bk [X]$ of regular functions on $X$
its $k$th infinitesimal neighborhood is the subscheme of $X$ with the defining ideal $I^k$ (in particular, an automorphism of such neighborhood
is essentially an automorphism of the algebra $\frac{A}{I^k}$).
It follows from \cite[Theorem 4.14 and Remark  4.16-4.17]{AFKKZ} that any automorphism $g: Z \to Z$ of an $k$th infinitesimal neighborhood $Z$
of  a zero dimensional subvariety  $Y$ of a smooth flexible variety $X$ is extendable to an automorphism of $X$
provided that $g|_Y =\id_Y$ and $g$ preserves local volume forms modulo $I^{k-1}$ (we say that such maps of infinitesimal neighborhoods
have Jacobian 1).
The second author  gave a proof  in \cite{Ud} of the similar result 
in the case when $Y$ is a smooth polynomial curve in a smooth flexible variety $X$ of $\dim X\geq 4$. 

In this paper we give a much stronger version of this result  (Theorem \ref{main.t1}) and develop some new important technical tools for the extension problem.
In particular, instead of polynomial curves we consider smooth curves in $X$ with trivial normal bundles.
(The class of such curves is sufficiently big.  It contains, for instance, curves which have only rational irreducible components.
 Furthermore, consider $X$ as a closed subvariety of an affine space.
Then by the Bertini theorem the intersection of $X$ with $n-1$ general hyperplanes where $\dim X=n$
is a smooth curve. Every irreducible component of this curve has a trivial normal bundle.)
As an application we get the following generalization of the Feller-van Santen theorem (Corollary \ref{stro.c1}).

\bthm\label{int.t1} Let $X$ be isomorphic to a  complex linear algebraic group  without nontrivial characters
and let $\dim X  \geq 4$.
Let $C_1$ and $C_2$ be smooth polynomial curves in $X$. 
Then every isomorphism of $k$th infinitesimal neighborhoods of $C_1$ and $C_2$ with Jacobian 1 extends to an automorphism of $X$.
\ethm

The absence of nontrivial characters is essential - without this assumption this theorem does not hold (see Remark \ref{stro.r1}).
If $C_1=C_2$, then  we can make a stronger statement in the holomorphic category without assuming that $C_i$ is
a polynomial curve or even that $C_i$ is irreducible. Namely, the holomorphic category has the advantage of
the Ivarsson-Kutzschebauch  theorem \cite{IK} which leads, in particular, to 
the following (Corollary \ref{hol.c2}).

\bthm\label{int.t2} Let $X$ be isomorphic to a  complex linear algebraic group  without nontrivial characters
and $\dim X  \geq 4$. Suppose that $C$ is a smooth closed curve in $X$ with a trivial normal bundle
and  defining ideal $\hI$ in the algebra $\Hol (X)$ of holomorphic functions on $X$. 
Let $\hA =\frac{\Hol (X)}{\hI}$ be the algebra of holomorphic functions on $C$.
Then  every $\hA$-automorphism $\theta : \frac{\Hol (X)}{\hI^k}\to \frac{\Hol (X)}{\hI^k}$   
 with Jacobian  1 extends to a holomorphic automorphism of $X$.

\ethm

In the case of $X\simeq \SL_n(\C)$ we can say even more (Corollary \ref{hol.c3}).

\bthm\label{int.t3} Let $X$ be isomorphic to $\SL_n (\C)$ where $n\geq 3$. Suppose that $C_1$ and $C_2$ are  isomorphic smooth closed curves in $X$
such that their normal bundles in $X$  are trivial.
Then every isomorphism of $k$th infinitesimal neighborhoods of $C_1$ and $C_2$ with Jacobian 1 extends to a
holomorphic  automorphism of $X$.
\ethm

In the case of sphere over $\bk$ the similar result holds in the algebraic category (Corollary \ref{stro.c2}).

\bthm\label{int.t4} Let $X$ be an $(m-1)$-sphere over $\bk$ where $m\geq 6$.  
Let $C_1$ and $C_2$ be isomorphic  smooth closed  curves in $X$ such that their irreducible components are rational.
Then every isomorphism of $k$th infinitesimal neighborhoods of $C_1$ and $C_2$ with Jacobian 1 extends
to an automorphism of $X$.\ethm

Let us briefly discuss the crucial step behind these theorems and the idea of its proof. As we mentioned, we
consider a more general situation when $X$ is a flexible variety of $n:=\dim X \geq 4$,  $C$ is a smooth closed 
curve in $X$ with a trivial normal bundle $N_XC$ and defining ideal $I\subset \bk[X]$. The consideration can be reduced to the
case of an irreducible $C$. Then $A=\bk [C]\simeq \frac{\bk[X]}{I}$ is a Dedekind domain and $N_XC$
can be viewed as a free $A$-module with a basis $v_1, \ldots , v_{n-1}$. 
In particular, any $A$-automorphism of $N_XC$ can be presented by an element of $\GL_{n-1} (A)$ where the latter group contains
a subgroup $\EE_{n-1} (A)$ generated by elementary matrices of the form $\id + a e_{ij}$ (where $\id$ is the identity matrix, $a\in A$, $1\leq i\ne j\leq n-1$ and $e_{ij}$
is the matrix with entry 1 in position $(i,j)$ and zeros everywhere else). 
The crucial step  is  that  under certain mild assumptions one has the following claim.

{\em  Every $A$-automorphism of $N_XC$ that is presented by an element of $\EE_{n-1} (A)$ is induced by an automorphism of $X$}.

Assume for simplicity that $X$ is smooth affine and that the natural action of a one-parameter unipotent subgroup of $\Aut (X)$ on $X$
(which is also called a $\G_a$-action) has a finitely generated ring $A^{\G_a}$ of invariants.
That is, there exists the algebraic quotient morphism $\rho : X\to Q=\Spec A^{\G_a}$ of this action sending $X$ to the normal affine variety $Q$. 
Recall that such a $\G_a$-action can be viewed as the flow of a so-called locally nilpotent vector field $\delta$ on $X$. 
Let $C'$ be the closure of $\rho (C)$ in $Q$ and $f$ be a regular function on $Q$ which vanishes on $C'$.
Then $\rho^*(f)\delta$ is again a locally nilpotent vector field whose flow $\exp (t\rho^*(f) \delta)$ leaves every point of $C$ fixed and, in particular,
it induces  $A$-automorphisms of $N_XC$. The plan is to show that every $A$-automorphism of $N_XC$ given by an
element of $\EE (A)$ is induced by a composition of  automorphisms of $X$ of the form $\exp (t\rho^*(f) \delta)$ with $\delta$ running over the set of locally
nilpotent vector fields on $X$.
To follow through this plan we use the fact (\cite[Theorem 6.1]{Ka20}) that  $\delta$ can be chosen so that for a given finite subset $K$ of $C$
there is an open subset $V'$ of $C'$ for which $V=\rho^{-1}(C')\cap C$ contains $K$ and $\rho|_V : V \to V'$ is an isomorphism.
Let $C\setminus V$ be the zero locus of an ideal $\fq\subset A$ in $C$, 
$\FF_{n-1}(A, \fq)$ be the subgroup of $\EE_{n-1}(A)$  generated by elementary matrices of the form $\id + a e_{ij}$ where $a \in \fq$
and $\EE_{n-1} (A, \fq)$ be the normal closure of $\FF_{n-1}(A, \fq)$ in $\EE_{n-1}(A)$.
Consider sections $w_1, \ldots, w_{n-1}$ of $N_XC$ such that they generate the free $\bk [V]$-module $N_XV$ and let $M$ be the 
free $A$-submodule  of $N_XC$ 
generated by these sections. Then  abundance of $\G_a$-action on $X$ enables us to prove that shrinking $V$ one can suppose that
for an appropriate choice of  $w_1, \ldots, w_{n-1}$  and sufficiently large $k>0$ every  $A$-automorphism of $M$ presented by a matrix from $\FF_{n-1} (A, \fq^{k})$
with respect to the basis  $w_1, \ldots, w_{n-1}$  is
induced by a composition of automorphisms of the form $\exp (t\rho^*(f) \delta)$ (see, Corollary \ref{main.c2}).
By \cite[Theorem 3]{Ni} $\FF_{n-1} (A, \fq^{k})\supset \EE_{n-1} (A, \fq^{2k})$.
Now one needs to switch from the basis $w_1, \ldots, w_{n-1}$ in $M$ to the basis $v_1, \ldots, v_{n-1}$ in $N_XC$
which  requires  some results about matrices over
Dedekind domains from  \cite{BMS} and \cite{Ba} (see, Theorem \ref{bass.t1}).
Under this change of bases the group $\EE_{n-1} (A, \fq^{2k})$ is transformed to a group containing
$\EE_{n-1} (A,\fq^{2mk})$ for some $m> 0$. Hence, 
every  $A$-automorphism of $N_XC$ presented by a matrix from $\EE_{n-1} (A,\fq^{2mk})$
with respect to the basis  $v_1, \ldots, v_{n-1}$  is
induced by a composition of automorphisms of the form $\exp (t\rho^*(f) \delta)$ (see, Corollary \ref{main.c2}).
Then, replacing $K$ by $C\setminus V$ in the argument above, we show that for some ideal $\fp$ in $A$ such that $\fq+\fp=A$
every  $A$-automorphism of $N_XC$ presented by a matrix from $\EE_{n-1} (A,\fp^{2mk})$ is again induced by an automorphism of $X$.
Observing, that $\EE_{n-1} (A,\fq^{2mk})$ and $ \EE_{n-1} (A,\fp^{2mk})$ generate $\EE_{n-1} (A)$ we get the desired claim.

With this claim at hand we note that for  an $A$-automorphism  $\theta: \frac{\bk [X]}{I^2} \to \frac{\bk [X]}{I^2}$ 
of the second infinitesimal neighborhood of $C$ the assumption that the Jacobian is 1 is
equivalent to the fact that the induced automorphism of $N_XC$ is contained in $\SL_{n-1}(A)$.
This implies that if the claim before is valid, then $\theta$ is induced by an automorphism of $X$ provided $\SL_{n-1}(A)=\EE_{n-1} (A)$. Say, rational curves satisfy this property
(which is the reason why such curves appear in Theorem \ref{int.t4}),
whereas in the holomorphic category one has always $\SL_{n-1}(\hA)=\EE_{n-1} (\hA)$ by virtue of the Ivarsson-Kutzschebauch  theorem
(see, Proposition \ref{hol.p1}).  To extend this result to automorphisms $\tilde \theta: \frac{\bk [X]}{I^k} \to \frac{\bk [X]}{I^k}$ of  the $k$th infinitesimal neighborhood of $C$
we use induction on $k$ with the case $k=2$ as the first step of induction. For implication $k-1 \implies k$
it suffices to show that $\tilde \theta$ is induced by an automorphism of $X$ if its restriction 
to the $(k-1)$st infinitesimal neighborhood of $C$ is trivial. The proof of the latter statement essentially copies  the similar step  in the proof of
 \cite[Theorem 4.14]{AFKKZ}.

The paper is organized as follows.  In Section 1 we collect technical tools developed in \cite{Ka20}
that are important for this paper.
Section 2 contains some criterion which enables us to deal with
curves that are not once-punctured (recall that a curve $C$ is once-punctured
if for any completion $\bC$ of $C$ the set $\bC\setminus C$ is a singleton). 
The advantage of a once-punctured curve (like $\C$) 
is that every morphism of this curve
to a quasi-affine algebraic variety is automatically proper.  Meanwhile, the assumptions under which we can prove the crucial step described before
include the properness of the morphisms $\rho|_C : C \to Q$.
The criterion we developed guarantees that for any closed curve $Z$  in a flexible variety $X$
and a morphism $\tau : X \to P$ to an affine variety $P$ there is an automorphism $\alpha$ of $X$ such that
$\tau\circ\alpha|_Z : Z \to Q$ is proper.
In Section 3  we fix some notations for the rest of the paper and spell out a convention under which our criterion of properness works.
In Section 4 using this convention we describe the structure of normal bundles of smooth curves in a large class of flexible varieties.
Actually, we show later (Corollary \ref{main2.c1}) that if a normal bundle of a smooth closed curve in such variety is trivial, then every nonvanishing 
section of this bundle is induced by a locally nilpotent vector  field (cf., \cite[Corollary 4.3]{AFKKZ}).
Section 5 contains some facts about special linear groups and their subgroups generated by elementary transformations over commutative rings. 
In Sections 6 we prove the main theorem (Theorem \ref{main.t1}) that states that
under  the assumption on Jacobians and some  other mild assumptions automorphisms of $k$th infinitesimal neighborhoods  of  smooth closed irreducible curves  in a flexible 
variety $X$ whose normal bundles are trivial  extend to automorphisms of $X$. In Section 7 we remove the assumption of irreducibility of these curves in
the main theorem.
As applications we prove Theorems \ref{int.t1} and \ref{int.t4} in Section 8. Section 9 contains Theorems \ref{int.t2} and \ref{int.t3}.

{\em Acknowledgements.} The authors are grateful to Leonid Vaserstein and Mikhail Zaidenberg for useful consultations.
They are also deeply indebted to the referee for a very thorough and helpful  review.
In particular, the referee corrected a mistake in the original proof of Lemma \ref{main.l7} and suggested a better proof
of Proposition \ref{hol.p1} which we use in this text.

\section{Preliminaries}

In this section we present some technical tools developed in \cite{Ka20} which require the following
definitions.

\bdefi\label{def.d1}  (1) Given an irreducible algebraic variety $\cA$ and
a map $\varphi:\cA\to\Aut(X)$  (where $\Aut (X)$ is the group of algebraic automorphisms of $X$) we say that $(\cA,\phi)$
is an {\em algebraic family of automorphisms of $X$} if the induced map
$\cA\times X\to X$, $(\alpha,x)\mapsto \varphi(\alpha).x$ is a morphism (see \cite{Ra}).

(2)  If we want to emphasize additionally that  $\varphi (\cA)$ is contained in a subgroup $G$ of $\Aut (X)$, then we say that
$\cA$ is an {\em algebraic  $G$-family} of automorphisms of $X$.            

(3) In the case when $\cA$ is a connected algebraic group and the induced map 
$\cA\times X\to X$ is not only a morphism but also an action of $\cA$ on $X$ we call this family a {\em connected algebraic subgroup} of $\Aut (X)$.
\edefi

\bdefi\label{def.d2}
We call a subgroup $G$ of $\Aut (X)$ {\em algebraically generated} if it is generated as an abstract group by a family 
$\cG$ of connected algebraic subgroups of $\Aut (X)$.
\edefi

\bdefi\label{def.d3}  (1) A nonzero derivation $\delta$ on the ring $A$ of regular functions on an affine algebraic variety $X$ is called {\em  locally nilpotent}
if for every $0\ne a \in A$ there exists a natural $n$ for which $\delta^n (a)=0$. 
This derivation can be viewed as a vector field on $X$ which
we also call {\em locally nilpotent}. The set of all locally nilpotent vector fields on $X$ will be denoted by $\LND (X)$. 
The flow of $\delta \in \LND (X)$ 
 is an algebraic $\G_a$-action on $X$, i.e., the action of the group $(\bk, +)$ 
which can be viewed as a one-parameter unipotent group $U$ in the group $\Aut (X)$ of all algebraic automorphisms of $X$.
In fact, every $\G_a$-action is a flow of a locally nilpotent vector field (e.g, see \cite[Proposition 1.28]{Fre}).

(2) Note that the kernel $\Ker \delta $ of $\delta$ coincides with the ring $\bk[X]^{\G_a}$ of $\G_a$-invariants in $\bk [X]$
where $\bk[X]^{\G_a}$ can be always treated as the ring of regular functions on a quasi-affine variety $V$ \cite[Theorem 1]{Wi1}).
The morphism $\tau : X \to V$ induced by the natural embedding 
$\bk[X]^{\G_a} \hookrightarrow \bk[X]$ is called the {\em categorical quotient morphism} of the $\G_a$-action.
When $V$ is an affine variety (i.e.,  when $\bk[X]^{\G_a}$ is finitely generated) 
we call $\tau$ an {\em algebraic quotient morphism}. 
Since $V$ is not always affine  by Nagata's example, we use sometimes  instead of $\tau$ partial quotient morphisms. 
Namely, a morphism  $\rho : X \to Q$ into a normal affine variety $Q$ is called a {\em partial quotient morphism associated with $\delta$} if
the ring $\rho^* (\bk [Q])$ is contained in $\Ker \delta$ and general fibers of $\rho$ are isomorphic
to $\A^1$.  If $X$ is normal, then partial quotient
morphisms always exist by virtue of the Rosenlicht Theorem (e.g., see \cite[Theorem 2.3]{PV} and \cite[Definition 2.16]{FKZ}). 

(3) If $X$ is a quasi-affine variety, then an algebraic vector field $\delta$ on $X$ is called {\em locally nilpotent} if $\delta$ extends
to a locally nilpotent vector field  $\tilde \delta$ on some affine algebraic variety $Y$ containing $X$ as an open subset such
that $\tilde \delta$ vanishes on $Y\setminus X$ where ${\rm codim}_Y (Y \setminus X) \geq 2$.  Note that under this assumption
$\delta$ generates a $\G_a$-action on $X$ and we use again the notation $\LND (X)$ for the set of all locally nilpotent vector
fields on $X$. As a partial quotient morphism $\rho : X \to Q$ of $\delta$ we consider the restriction of a partial quotient morphism of $\tilde \delta$ to $X$.
\edefi

\brem\label{def.r1} Definition \ref{def.d3} (2)-(3) implies that  
for a  dense open smooth subvariety 
$Q_0$ of $Q$  and $X_0=\rho^{-1}(Q_0)$  the morphism $\rho|_{X_0} : X_0\to Q_0$ is smooth and surjectve
(in fact, one can assume that $X_0$ is naturally isomorphic to $Q_0\times \A^1$  by \cite[Theorem 2]{KaMi})
and  $\delta$ does not vanish on $X_0$. 
\erem

\bdefi\label{def.d4}

(1) For every locally nilpotent vector fields $\delta$ and each function $f \in \Ker \delta$ from its kernel the field
$f\delta$ is called a replica of $\delta$. Recall that such replica is automatically locally nilpotent.

(2)  Let $\cN$ be a set of locally nilpotent vector fields on $X$ and $G_\cN \subset  \Aut (X)$ denotes the group generated
by all flows of elements of $\cN$. We say that $G_\cN$ {\em is generated by $\cN$}.

(3) A collection of locally nilpotent vector fields $\cN$ is called saturated if $\cN$ is closed under conjugation by elements in $G_\cN$
and for every $\delta \in \cN$ each replica of $\delta$ is
also contained in $\cN$.

\edefi

\bdefi\label{def.d5}  Let $X$ be a normal quasi-affine algebraic variety of dimension at least 2,
$\cN$ be a saturated set of locally nilpotent vector fields on $X$ and  $G=G_\cN$ be the group
generated by  $\cN$.
Then $X$ is called $G$-flexible if for any point $x$ in the smooth part $X_{\rm reg}$ of $X$ the vector space $T_xX$ is generated 
by  the values of locally nilpotent vector fields from $\cN$ at $x$
(which is equivalent to the fact that $G$ acts transitively on $X_{\rm reg}$ \cite[Theorem 2.12]{FKZ}). In the case of $G=\SAut (X)$ we call $X$ flexible
without referring to $\SAut (X)$ (recall that $\SAut (X)$ is the subgroup of $\Aut X$ generated by all one-parameter unipotent subgroups).
\edefi



The following Collective Transversality Theorem is one of the main technical tools in \cite{Ka20}.

\bthm\label{pre.t1} {\rm (\cite[Theorem 1.4]{Ka20}, cf. \cite[Theorem 1.15 ]{AFKKZ})}  Let $X$ and $P$ be smooth irreducible
algebraic varieties and $\kappa : X \to P$ be a smooth morphism.
Let $G$ be a subgroup of the group  $\Aut(X/P)$  of automorphisms of $X$ over $P$ such that $G$ is
algebraically generated by a system $\cG$ of connected algebraic
subgroups closed under conjugation in $G$. 
Suppose that the restriction of the $G$-action to $\kappa^{-1} (p)$ is transitive for every $ p\in P$.

Then there exist subgroups $H_1,\ldots, H_m\in \cG$  such that for
any locally closed reduced subschemes $Y$ and $Z$ in $X$ one can
find a  dense open subset $U=U(Y,Z)\subseteq H_1\times
\ldots \times H_m$ so that every element $(h_1,\ldots, h_m)\in
U$ satisfies the following:  

{\rm (i)} 
 $\dim (Y\cap (h_1\cdot\ldots\cdot h_m).Z)\le \dim (Y\times_PZ) + \dim P -\dim X$. 
 
 In particular, when  $\dim Y\times_PZ \leq \dim Y+ \dim Z -\dim P$ one has 

{\rm (ii)} $\dim (Y\cap (h_1\cdot\ldots\cdot h_m).Z)\le \dim Y+\dim Z -\dim X$. 

Furthermore, suppose that the inequality $\dim Y\times_PZ \leq \dim Y+ \dim Z -\dim P$ holds
and also that $Z$, $Y \times_PZ$, and  $Y\times_PX$ are smooth.
Then

{\rm (iii)}  $(h_1\cdot\ldots\cdot h_m).Z$ meets $Y$ transversally.

\ethm

The crucial step in the proof of Theorem \ref{pre.t1} is the follwoing.

\bprop\label{pre.p1} {\rm (\cite[Proposition 1.7]{Ka20})}
 Let the assumption of Theorem \ref{pre.t1} hold.
Then there is a sequence $H_1,\ldots, H_m$ in $\cG$
so that for a suitable open dense subset $U\subseteq
H_{m}\times\ldots \times H_{1}$, the map
\be\label{pre.eq1}
\Phi : H_m\times \ldots\times H_1\times X \lto X\times_P X
\quad\mbox{with} \quad (h_m,\ldots,h_1,x)\mapsto
((h_m\cdot\ldots\cdot h_1).x ,x)
\ee
is smooth on $U\times X$.

\eprop

Namely, we have the following.

\bprop\label{pre.p2}{\rm (\cite[Proposition 1.10]{Ka20})}
If a sequence $H_1,\ldots, H_m \in \cG$ satisfies the conclusions of Proposition \ref{pre.p1},  then it satisfies also the conclusions of
Theorem \ref{pre.t1}. Furthermore, for any element $H$ of $\cG$
the sequence  $H_1,\ldots, H_m, H$ (resp. $H, H_1,\ldots, H_m$) satisfies  the conclusions of Proposition \ref{pre.p1} and Theorem \ref{pre.t1} as well.
\eprop

We shall need the following notions which, unfortunately, were not introduced in \cite{Ka20}.

\bdefi\label{pre.d6} Let $\kappa : X \to P$, $G$ and $\cG$ satisfy the assumptions of Theorem \ref{pre.t1} and $\dim X=n$. 
Consider $(X\times_PX)\setminus \Delta$ (where $\Delta$ is the diagonal), the complement $T'X$ to the zero section in the tangent bundle
of $X$ and the frame bundle ${\rm Fr} (X)$ bundle of $TX$ (i.e., the fiber of  ${\rm Fr} (X)$ over $x \in X$ consists of all bases of $T_xX$).
Projectivization of $TX$ replaces ${\rm Fr}(X)$ with a bundle  ${\rm PFr} (X)$ whose fiber over $x$ consists of all ordered $n$-tuples of points
in the proectivization $\PP^n$ of $T_xX$ that are not contained in the same hyperplane of $\PP^n$. 
Then we have natural $G$-actions on all these objects.
Let $Y$ be either $X$, or $(X\times_PX)\setminus \Delta$, or $T'X$,  or ${\rm PFr} (X)$. Suppose that  the $G$-action is transitive on every fiber of $Y$ over $P$.
Then we say that an algebraic $G$-family $\cA$ of automorphisms of $X$ is
{\em a regular $G$-family  for $Y$ over $P$} if

(i) $\cA= H_{m}\times\ldots \times H_{1}$ where each $H_i$ belongs to $\cG$;

(ii) for a suitable open dense subset $U\subseteq H_{m}\times\ldots \times H_{1}$, the map
\be\label{pre.eq2}  \Psi: H_m\times \ldots\times H_1\times Y \lto Y\times_P Y
\quad\mbox{with} \quad (h_m,\ldots,h_1,y)\mapsto
((h_m\cdot\ldots\cdot h_1).y ,y)
\ee
is smooth on $U\times Y$. 

Of course, by Proposition \ref{pre.p1} the assumption about transitivity implies that regular families exist.
Note also that if $\cA$ is a regular $G$-family, say, for $(X\times_PX)\setminus \Delta$ over $P$ and $\cB$ is a regular $G$-family for $T'X$ over $P$,
then by Proposition \ref{pre.p2} $\cA \times \cB$ (resp. $\cB\times \cA$) is a regular $G$-family for both $(X\times_PX)\setminus \Delta$ and $T'X$
over $P$. This implies the existence of algebraic $G$-families $\cA$ that are regular for all four varieties
 $X$,  $(X\times_PX)\setminus \Delta$,  $T'X$ and ${\rm PFr} (X)$ over $P$. We call such $\cA$ a {\em perfect $G$-family for $Y$ over $P$}
 (if $P$ is a singleton, then we just say that $\cA$ is a {\em  perfect $G$-family for $Y$}).

\edefi

\brem\label{pre.r1} It is worth observing the following property of regular $G$-families in the case when $P$ is a singleton. Let $K$ be a finite subset of $Y$,
$Z$ be a dense open subset of $Y\times K$ and $Z_y =Z\cap (Y\times y)$ where $y \in K$.
Note that $V_y=(\Psi|_{U\times y})^{-1}(Z_y)$ is  a nonempty open 
subset of $U$ since $\Psi$ is smooth. Hence,  $\Psi(\bigcap_{y \in K}V_y\times K)$ is a nonempty open 
subset of $Y\times K$.
This implies that the morphism $$U \to Y\times K, \, \alpha \mapsto (\alpha (y))_{y\in K}$$ is dominant.

\erem

\bexa\label{pre.exa1} Let $X$ be a smooth $G$-flexible variety where $G$ 
is generated by a saturated set $\cN \subset \LND (X)$.  
Then the $G$-action is $m$-transitive for every $m>0$ \cite[Theorem 2.12]{FKZ}.
Hence, the  natural action of $G$ on $(X\times_PX)\setminus \Delta$ is transitive.  Similarly, 
treating $X$ as an open subset of a normal affine variety one can apply \cite[Theorem 4.14 and Remark 4.16]{AFKKZ} to show that
$G$ acts transitively on $T'X$ and on ${\rm PFr} (X)$ (but not on ${\rm Fr} (X)$). 
Hence, $X$ admits a perfect $G$-family.
\eexa

Theorem \ref{pre.t1} was used to establish the next fact.  

\bthm\label{pre.t2}  {(\cite[Theorem 4.2(i)(iii)]{Ka20}) }
Let $X$ and $P$ be smooth algebraic varieties and $Q$ be a normal algebraic variety.  Let $\rho : X \to Q$ and $\tau: Q\to P$ be dominant morphisms such that $\kappa : X \to P$ is smooth for $\kappa=\tau \circ \rho$. 
Suppose that  $Q_0$ is a non-empty  open smooth subset of $Q$ so that for  $X_0=\rho^{-1} (Q_0)$ the morphism $\rho|_{X_0} : X_0 \to Q_0$ is smooth.
Let $G$ be an algebraically generated subgroup of the group  $\Aut(X/P)$  of automorphisms of $X$ over $P$ such that $G$ acts
2-transitively on each fiber of $\kappa : X \to P$ and let
$Z$ be a locally closed reduced subvariety in $X$. 

{\rm (i)} Let $\dim Z\times_P Z\leq 2\dim Z- \dim P$ and $\dim Q \geq \dim Z+m$ where $m \geq 1$. Then  there exists an algebraic $G$-family $\cA$ of automorphisms
of $X$  such that
for a general element $\alpha \in \cA$ one can find 
a constructible subset $R$ of $\alpha (Z)\cap X_0$ of dimension $\dim R \leq \dim Z-m$ for which
$\rho (R)$ and $ \rho (\alpha (Z)\setminus R)$ are disjoint  
and the restriction $\rho|_{(\alpha (Z)\cap X_0)\setminus R}:  (\alpha (Z) \cap X_0) \setminus R\to  Q_0$ 
of $\rho$ is injective.

In particular, if $\dim Q\geq  2\dim Z +1$  and $\bZ_\alpha'$ is the closure of  $Z_\alpha'=\rho \circ \alpha (Z)$ in $Q$, then
for a general element $\alpha \in \cA$ the map $\rho|_{\alpha (Z)\cap X_0}:  \alpha (Z) \cap X_0 \to Z_\alpha'\cap Q_0$ is a bijection,
while in the case of a pure-dimensional $Z$ and $\dim Q \geq \dim Z +1$ 
the morphism $\rho|_{\alpha (Z)\cap X_0}:  \alpha (Z) \cap X_0 \to \bZ_\alpha'\cap Q_0$ is birational.

{\rm (ii)} Suppose that $X$ is $G$-flexible where
$G$ is generated by a saturated set $\cN\subset \LND (X)$ and $P_0 = \tau (Q_0)$. 
Let $\dim Z\times_P Z\leq 2\dim Z- \dim P, 2\dim Z +1 \leq \dim Q$ and $\dim T(Z_0/P_0)  \leq \dim Q -\dim P$.
Then the algebraic family $\cA$ from (i) can be chosen so that 
for a general element $\alpha \in \cA$   the morphism  $\rho|_{\alpha (Z)\cap X_0}:  \alpha (Z) \cap X_0 \to Z_\alpha'\cap Q_0$ is 
bijective and it induces an injective map of the  tangent bundle of $\alpha (Z) \cap X_0$ into the  tangent bundle of  $Q_0$.

\ethm

\brem\label{pre.r2}  (1) We shall need below the case when $P$ is a singleton. In particular, the assumption $\dim Z\times_P Z\leq 2\dim Z- \dim P$ 
can be omitted.

(2) If $X$ is a smooth $G$-flexible variety, then the $G$-action on $X$ is $m$-transitive for every $m\geq 1$ \cite[Theorem 2.12]{FKZ}.
In particular, when $P$ is a singleton the assumption that  $G$ acts 2-transitively on $X$ can be also omitted in this case.

(3)  It follows from the proof  that  any regular $G$-family for $(X\times_PX)\setminus \Delta$ over $P$ can serve as $\cA$  
in Theorem \ref{pre.t2} (i), whereas for Theorem \ref{pre.t2} (ii) any  $G$-family regular for both $(X\times_PX)\setminus \Delta$ and $T'X$ over $P$
works. Thus, whenever we apply Theorem \ref{pre.t2} below we use a perfect $G$-family.
In particular,  for every $H \in \cG$ the families $H\times \cA$ and $\cA \times H$ also satisfy
the conclusions of Theorem \ref{pre.t2}.  
\erem

Another fact from \cite{Ka20} (also based on Theorem \ref{pre.t1}) which we shall need later is the following.

\bthm\label{pre.t4}{\rm (\cite[Theorem 6.1]{Ka20})} Let $X$ be a smooth quasi-affine algebraic variety, $\cN$ be a saturated set of locally nilpotent vector fields on $X$, 
and $G\subset \SAut (X)$ be the group generated by $\cN$. Suppose that $X$ is $G$-flexible.  
Let $\rho  : X \to Q$ be a partial quotient morphism
associated with a nontrivial $\delta \in \cN$, $Z$ be a
locally closed reduced subvariety of $X$
of codimension at least 2 and $K$ be a finite subset of $Z$ such that $\dim T_{z_0}Z \leq \dim Q$ 
for every $z_0 \in K$.
Then there exists a connected algebraic $G$-family  $\cA$ of automorphisms of $X$
such that for a general element $\alpha \in \cA$
and  the closure $\bZ_\alpha'$ of $Z_\alpha'=\rho\circ \alpha (Z)$ in $Q$ 
one can find a neighborhood $V_{\alpha}'$ of $\rho (\alpha (K))$ in $\bZ_\alpha'$ such that for $V_{\alpha}=\rho^{-1}(V_{\alpha}')\cap \alpha (Z)$ 
the morphism $\rho|_{V_{\alpha}} : V_{\alpha} \to V_{\alpha}'$ 
is an isomorphism.

\ethm

\brem\label{pre.r3} (1) The proof of Theorem \ref{pre.t4} (see \cite[page 550]{Ka20})  implies that $\cA$ is of the form $\cA =U_\kappa\times \cA_1$
where $U_\kappa$ is a fixed family\footnote{This family $U_\kappa$ is a product of  several one-parameter unipotent subgroup depending 
on $\delta$ and $K$ but not on $Z$ \cite[Formula (20)]{Ka20}.}, 
while $\cA_1=H_1\times \ldots \times H_m$ is any regular $G$-family for $X$ that contains an element $\alpha_0\in G$
with the following property:  for every $z\in \alpha_0(K)$ and $q=\rho (z)$ the induced map $\rho_* : T_z\alpha (Z) \to T_qQ$
is injective. Note that by Remark \ref{pre.r1} such $\alpha_0$ exists automatically if one requires that $\cA_1$ is  a regular $G$-family
for ${\rm PFr} (X)$. Thus, from now on whenever we apply Theorem \ref{pre.t4} we suppose that $\cA_1$ is a perfect family.

(2) Another important feature of Theorem \ref{pre.t4} we want to emphasize is that by construction 
$V_\alpha'$ is contained in $Q_0\subset Q$ as in Remark \ref{def.r1}. 

\erem

In conclusion of this section we state the following simple fact which will be crucial in further considerations. 

\blem\label{pre.l1} {\rm(\cite[Lemma 4.1]{AFKKZ})}\footnote{In \cite{AFKKZ} this fact was formulated for an affine variety $X$ but the proof goes without
change in the quasi-affine case.}
Let $\delta$ be a locally nilpotent vector field on a smooth quasi-affine variety $X$, and let  $p\in X$ be a point. 
Assume that $f\in \Ker \delta$  and $f(p)=0$. If  $\Phi=\exp (f\delta )$ is the automorphism associated with the replica $f\delta$, then
\be\label{pre.eq3}   \dd_p \Phi (w)=w+\dd_p f(w) \delta (p)\,\,\, \text{for all} \,\,\, w\in T_pX      \ee
where $\delta (p) \in T_pX$ is the value of the vector field $\delta$ at $p$.
\elem

\section{Properness}

\bdefi\label{prop.d1}  Let $\delta_1$ and $\delta_2$ be locally nilpotent vector fields on an affine algebraic variety $X$ 
and $\rho_i : X\to Q_i$ be partial quotient morphisms associated with $\delta_i, \, i=1,2$. 

{\rm (1)} The pair $\delta_1, \delta_2$ will be called {\em suitable}  
if $\rho_1$ and $\rho_2$ can be chosen so that the morphism $\tau =(\rho_1, \rho_2) : X \to Q_1\times Q_2$ 
is proper. 

{\rm (2)}  The pair $\delta_1, \delta_2$ will be called  {\em semi-compatible}  if $\Ker \delta_1 $ and $\Ker \delta_2 \subset \bk [X]$ are finitely generated
(i.e., $\rho_1$ and $\rho_2$ can be chosen to be algebraic quotient morphisms), $\tau  : X \to \tau (X)$  is  birational and finite where $\tau (X)$ is closed in $ Q_1\times Q_2$. 
In particular, every semi-compatible pair of locally nilpotent vector fields is automatically suitable.
\edefi

It is useful to keep in mind the following criterion of semi-compatibility.

\bprop\label{prop.p1} {\rm (\cite[Proposition 3.4]{KaKu08})}  Let $\delta_1$ and  $\delta_2$ be locally nilpotent vector fields on a normal affine algebraic variety $X$
such that their kernels $\Ker \delta_1 $ and $\Ker \delta_2 \subset \bk [X]$ are finitely generated.  
Then these vector fields are semi-compatible if and only if  the span of $\Ker \delta_1 \cdot \Ker \delta_2$ contains a nontrivial ideal of $\bk [X]$.
\eprop

Though the next fact is not be used later in the paper, we include it as a natural property of suitability.

\bprop\label{prop.p2}  {\rm (cf.  \cite[Proposition 3.6]{KaKu08})}
Let $\delta_1$ and $ \delta_2$ be suitable locally nilpotent vector fields on  a affine algebraic variety $X$
with finitely generated kernels $\Ker \delta_1 $ and $\Ker \delta_2$.  
Suppose that $\Gamma$ is a finite group acting on $X$ so that each
$\delta_i$ induces a locally nilpotent vector field $\delta_i'$ on the quotient $X'=X/\Gamma$.
Then $\delta_1'$ and $\delta_2'$ are suitable on $X'$.
\eprop

\bproof Let $\rho_i: X \to Q_i$ and $\tau$ be as in Definition \ref{prop.d1} and let $K_i\subset \SAut (X)$ be the flow of $\delta_i$.
Since by the assumption the  $K_i$-action on $X$ commutes  with the $\Gamma$-action,  we have a $\Gamma$-action on $Q_i$
such that $\rho_i$ is $\Gamma$-equivariant.
We have also the induced morphisms $\rho_i': X' \to Q_i':=Q_i/\Gamma$. 
Every morphism $\varphi' : X' \to Z$ that is constant on the orbits of $K_i$ induces
a $\Gamma$-equivariant morphism $\varphi : X \to Z$ (where $Z$ is equipped with a trivial $\Gamma$-action). Since $\rho_i$
is an algebraic quotient morphism, $\varphi $ is the composition of $\rho_i$ and a morphism $\psi : Q_i \to Z$. Since $\rho_i$ is $\Gamma$-equivariant,
so is $\psi$. Hence, $\psi$ factors through $Q_i'$ which shows that $\rho_i' : X'\to Q_i'$ is the algebraic quotient morphism. 
Consider the following commutative diagram

\[  \begin{array}{ccc} X &  \stackrel{{ \tau}}{\rightarrow} & Q_1\times Q_2 \\
\, \, \, \, \downarrow {\kappa}  & & \, \, \, \, \downarrow { \kappa'} \\
X' &  \stackrel{{ \tau'}}{\rightarrow} & Q_1'\times Q_2'

\end{array} \]

\noindent where $\tau' =(\rho_1',\rho_2')$, $\kappa$ is the quotient morphism of the $\Gamma$-action on $X$ and
$\kappa'$ is the quotient morphism of the natural $\Gamma \times \Gamma$-action on $Q_1\times Q_2$.
Since $\tau$ and $\kappa'$ are proper morphisms, so is $\kappa' \circ \tau=\tau'\circ \kappa$. Since $\kappa$ is surjective, this
implies that $\tau'$ is proper which is the desired conclusion.
\eproof

The next simple observation is crucial.

\bthm\label{prop.t1}  Let  $\{ \delta_\beta | \, \beta \in \cB\}$ (where $\cB$ is an index set) be an infinite collection of
locally nilpotent vector fields on an affine algebraic variety $X$ such that for every $\beta \ne \gamma \in \cB$  the pair $(\delta_\beta, \delta_\gamma)$
is suitable for appropriate partial quotient morphisms $\rho_\beta : X \to Q_\beta$ and $\rho_\gamma : X \to Q_\gamma$. 
Suppose that 
$C$ is a closed (but not necessarily irreducible) curve in $X$. Then there is a finite subset $T$ of
$\cB$ such that for every $\beta \in \cB \setminus T$ the morphism $\rho_\beta|_C : C \to Q_\beta$ is proper.
\ethm

\bproof 
Suppose that $\bC$ is a completion of $C$ such that it is smooth at every point of $\bC\setminus C$.
By \cite[Chapter II, Theorem 4.7]{Har} the fact that a morphism $\theta : C \to Y$ is not proper is equivalent to the fact that $\theta$ extends regularly to 
some point $c \in \bC \setminus C$.  Assume that $\rho_\beta|_C : C \to Q_\beta$ extends to $c \in \bC \setminus C$.
Then for every $\beta \ne \gamma \in  \cB$ the morphism $\rho_\gamma|_C : C \to Q_\gamma$ cannot be extended to $c$.
Indeed, otherwise the morphism $\tau|_C : C \to Q_\beta\times Q_\gamma$ extends to $c$ where $\tau=(\rho_\beta, \rho_\gamma) : X \to Q_\beta\times Q_\gamma$.
That is, $\tau$ is not proper by \cite[Chapter II, Corollary 4.8]{Har} which contradicts suitability.
Thus, there is at most a finite number of elements $\beta$ of $ \cB$ such $\rho_\beta|_C$ extends to some $c \in \bC\setminus C$.
This yields the desired conclusion.
\eproof

\bdefi\label{prop.d2}
Let $X$ be an affine algebraic variety.
We say that $X$ is {\em parametrically suitable} if there exists a connected algebraic family 
$\cB$ of automorphisms of $X$ with $\dim \cB \geq 1$, a locally nilpotent vector field 
$\delta \in \LND (X)$ and a collection $\{ \rho_\beta : X \to Q_\beta | \, \beta \in \cB \}$ of partial quotient morphisms $\rho_\beta$
associated with  $\delta_\beta: =\beta_* (\delta)$  such that for every $\beta\ne \gamma \in \cB$ 
the morphism $(\rho_\beta, \rho_\gamma) : X \to Q_\beta \times Q_\gamma$ is proper (i.e.,
the  pair $(\delta_\beta, \delta_\gamma)$ is suitable).
Suppose additionally that $G \subset \Aut (X)$ is a group generated by a saturated set $\cN \subset \LND (X)$. 
If $\delta$ belongs to $\cN$ and $\cB$ is contained in $G$, then we say that $X$ is {\em parametrically $G$-suitable}.
If we want to emphasize that we use a specific $\delta$ in this definition we speak about  parametric $G$-suitability  {\em with respect to $\delta$}.
\edefi

\brem\label{prop.r1}

(1) For $\delta \in \LND (X)$, $\alpha \in \Aut (X)$ and  $\delta_\alpha =\alpha_* (\delta)$ one has by definition  $\delta_\alpha (f )=  \delta (f\circ \alpha)$ for every $f \in \bk [X]$.
Hence, if $g=f\circ \alpha\in \Ker \delta$, then  $f=g\circ \alpha^{-1}\in \Ker \delta_\alpha$.  This implies that if 
$\rho : X \to Q$ is a partial quotient morphism associated with $\delta$, then $\rho\circ \alpha^{-1} : X \to Q$
is a partial quotient morphism associated with $\delta_\alpha$. In particular, we shall look for elements of a collection
$\{ \rho_\beta : X \to Q_\beta | \, \beta \in \cB \}$ in Definition \ref{prop.d2} in the form $\rho_\beta =\rho \circ \beta^{-1}$ (i.e., $Q_\beta =Q$).

(2)  Note also that one can suppose that $\cB$ in Definition \ref{prop.d2} contains the identity automorphism.
Indeed, if a pair $(\delta_\beta, \delta_\gamma )$ is suitable, then the pair $(\delta_{\alpha \circ \beta}, \delta_{\alpha \circ \gamma} )$
is also suitable, since $(\rho_{\alpha \circ \beta}, \rho_{\alpha \circ \gamma}) =(\rho_\beta, \rho_\gamma)\circ \alpha^{-1}$. 
Hence, $\cB$ can be replaced by $\alpha \cB$. Choosing $\alpha =\beta^{-1}$ for some $\beta \in \cB$ we get the desired statement.

(3) One can modify the notion of suitability and consider it not for pairs but, say, for triples of
locally nilpotent vector fields on $X$. That is, 
for a connected nonconstant algebraic family $\cB$ of automorphism of $X$, $\delta \in \LND (X)$,
a collection $\{ \rho_\beta : X \to Q_\beta | \, \beta \in \cB \}$ of partial quotient morphisms $\rho_\beta$
associated with  $\delta_\beta$
and every triple $(\beta, \gamma, \chi)$ of distinct elements of  $\cB$ one can require that the morphism $(\rho_\beta, \rho_\gamma,\rho_\chi): X \to Q\times Q\times Q$
is proper.  Then the straightforward adjustment of the proof of Theorem \ref{prop.t1} shows that $\rho_\beta|_C : C \to Q$ is proper for
a general $\beta \in \cB$. In this case we call $X$ {\em weakly parametrically suitable} (we do not know if this notion is equivalent to
parametric suitability).

(4) Though we use partial quotient morphisms in the definition of (parametric) suitability in all applications below
$\delta$ (and, consequently, $\delta_\beta$) admits an algebraic quotient morphism.

\erem

The following is one of the main criteria for parametric suitability.

\bthm\label{prop.t2} Let $X$ be a smooth complex affine algebraic variety, $F\simeq \SL_2(\C)$ be  an algebraic subgroup of $\Aut (X)$ contained in
a group $G\subset \Aut (X)$ such that the natural $F$-action on $X$ is fixed point free and
non-degenerate\footnote{An action of a linear algebraic group is non-degenerate if general orbits have the same dimension as the group.}. 
Then $X$ is parametrically $G$-suitable.
\ethm

\bproof  Fixing some basis $\bar u=(u_1,u_2)$ in $\C^2$ with the natural action of $\SL_2(\C)$ we get 
a unipotent subgroup $H_1$ (resp. $H_2$) of upper (resp. lower) triangular matrices in $\SL_2(\C)\simeq F$. 
The $\G_a$-action induced by $H_1$ (resp. $H_2$) on $X$ is
the flow of a locally nilpotent vector field $\delta_1$ (resp. $\delta_2$).
 By \cite[Theorem 3.1]{Had} there exists an algebraic quotient morphism $\rho : X \to Q$ associated with $\delta_1$.
The central fact we use now is that
$\delta_1$ and $\delta_2$ are semi-compatible and even compatible in the terminology of \cite[Theorem 12]{DDK}.
However, $H_1$ and $H_2$ (resp. $\delta_1$ and $\delta_2$) depend on the choice of $\bar u$.
In particular, for every $\alpha \in \SL_2(\C)$ the groups $\alpha H_1\alpha^{-1}$ and  $\alpha H_2\alpha^{-1}$
present upper and lower triangular matrices in an appropriate basis (note that this conjugation sends $\delta_i$ to $\alpha_*(\delta_i)$). 
Consider $\cB=H_2$ and let us show that for general $\alpha\ne \beta \in \cB$ the locally nilpotent vector fields
$\alpha_*(\delta_1)$ and $\beta_*(\delta_1)$ are semi-compatible, or, equivalently, 
in an appropriate basis $\alpha H_1\alpha^{-1}$ and $\beta H_1\beta^{-1}$ correspond to upper and lower unipotent matrices respectively.
The latter is equivalent to the similar claim about the groups $\beta^{-1} \alpha H_1\alpha^{-1}\beta$ and $H_1$.
We can suppose that $\alpha: (u_1,u_2)\mapsto (u_1, u_2+au_1)$ and  $\beta : (u_1,u_2)\mapsto (u_1, u_2+bu_1)$.
Then a direct computation shows that $\beta^{-1}\alpha H_1\alpha^{-1}\beta$ consists of matrices of the form
\[  \left[ {\begin{array}{cc}
 1-(a-b)t & t  \\      
-(a-b)^2t & 1+ (a-b)t  
 \end{array} } \right].  \]
 Another direct computation shows that for the automorphism $\gamma: (u_1,u_2) \mapsto (u_1 +(a-b)^{-1}u_2, u_2)$
 one has $\gamma H_2\gamma^{-1}= \beta^{-1}\alpha H_1\alpha^{-1}\beta$.  Hence, $\beta^{-1}\alpha H_1\alpha^{-1}\beta$ and $H_1$
 correspond to upper and lower unipotent matrices in an appropriate basis since $\gamma H_2\gamma^{-1}$ and $H_1=\gamma H_1\gamma^{-1}$ do.
 Thus,
 $\alpha_* (\delta_1)$ and $\beta_*(\delta_1)$ are suitable for general $\alpha\ne \beta \in \cB$ which yields the desired conclusion.
\eproof

\bexa\label{prop.ex1} Let $H$ be a complex semi-simple Lie group with Lie algebra different from $\fs\fl_2$.
Let $R$ be a proper closed reductive subgroup of $H$.  Then the homogeneous space $H/R$ is 
parametrically suitable. Indeed, by \cite[Theorem 24 and the proof of Corollary 25]{DDK} $H/R$ admits a non-degenerate fixed point free $\SL_2(\C)$-action.
Hence, we are done by Theorem \ref{prop.t2}.
\eexa

\bcor\label{prop.c2}
 Let $X$ be a smooth affine algebraic variety equipped with a fixed point free 
non-degenerate $\SL_2(\bk)$-action.  Let $\rho : X \to Q$ and $\cB$ be the same as in the proof of Theorem \ref{prop.t2} (but with $\C$ replaced by $\bk$).
Then for every closed curve $C$ in $X$ and  a general $\beta \in \cB$ the morphism $\rho_\beta|_C=\rho \circ \beta^{-1}|_C : C \to Q$ is proper.
\ecor

\bproof   If $\bk$ is a universal domain (i.e., it is an algebraically closed field of characteristic 0 with infinite transcendence degree over $\Q$),
then one can just refer to \cite{Ek}. To include a finite transcendence degree in the proof one can argue like this. 
Let $\tC$ be an open subset in a completion of $C$ consisting of $C$ and one extra point $z$. 
Consider $\tau_\beta : \tC\times \cB \dashrightarrow Q$ that sends each $(c,\beta) \in C\times \cB$ to $\rho_\beta (c)$. 
Assume that $\rho_\beta : C \to Q$ is not proper for
a general $\beta \in \cB$.  Then for some $\tC$ as above $\tau_\beta$ is a morphism.
By the ``finite extension principle" there exist an algebraically closed subfield $\bk^0$ of $\bk$ which is of a finite transcendence degree over of $\Q$, 
algebraic varieties $X^0,Q^0,\tC^0=C^0\cup z^0,\cB^0$ over $\bk^0$,
morphisms $\rho_\beta^0: C^0\to Q^0$,  $\tau_\beta^0: \tC^0\times \cB^0 \rightarrow Q^0$ and an $\SL_2(\bk_0)$-action on $X^0$ such that
$\rho_\beta: C\to Q$, $\tau_\beta : \tC\times \cB \to Q$ and the $\SL_2(\bk)$-action on $X$ are obtained by the base extension $\otimes_{\Spec \bk_0}\Spec \bk$.
By the Lefschetz principle we can view $\bk^0$ as a subfield of $\C$. The extension $\otimes_{\Spec \bk_0}\Spec \C$
leads to morphisms $\rho_\beta^\C: C^\C\to Q^\C$, $\tau_\beta^\C: \tC^\C\times \cB^\C \rightarrow Q^\C$ of complex varieties  and a non-degenerate 
fixed point free $\SL_2(\C)$-action on $X^\C$. However, Theorems \ref{prop.t1} and \ref{prop.t2} imply that $\rho_\beta^\C: C^\C\to Q^\C$
is proper for general $\beta \in \cB$ which contradicts the existence of the morphism $\tau_\beta^\C$. Hence, we have the desired conclusion.
\eproof

\section{Notations and Convention}

Let us list some notations which will remain unchanged through the rest of the paper.

\bsett\label{nor.s1} We suppose further that

- $X$ is a normal quasi-affine variety of dimension $n$;

- $\delta$ is a locally nilpotent nonzero vector field on $X$;

- $\rho : X \to Q$ is a partial quotient morphism associated with $\delta$;

- $Q_0$ is a smooth open dense subset of $Q$ such that $\rho|_{\rho^{-1} (Q_0)} : \rho^{-1} (Q_0)\to Q_0$ is a smooth morphism
and $\delta$ does not vanish on $\rho^{-1} (Q_0)$ (i.e., $Q_0$ is as in Remark \ref{def.r1});

- for any $\alpha \in \Aut (X)$  we let $\delta_\alpha =\alpha_* (\delta)$ and $\rho_\alpha =\rho \circ \alpha^{-1}: X\to Q$
(i.e., $\rho_\alpha$ is a partial quotient morphism associated with $\delta_\alpha$);

- $C$ is  a smooth  (but not necessarily irreducible) curve  in $X$ such that $C\subset X_{\rm reg}$ and $C$ is closed
in an affine variety containing $X$ as a dense open subset;

- $I$ is the defining ideal of $C$ in $\bk [X]$ (except for Section 5 where $I$ is an identity matrix);

- $\pr : TX|_C \to N_XC$ is the natural projection  onto the normal bundle $N_XC$ of $C$ in $X$.

\esett

In the case when $C$ contains an irreducible component which is not a once-punctured curve we shall often use
the following important convention.

\bconv\label{nor.conv1} Suppose additionally that $X$ is a $G$-flexible variety 
 with $\dim X=n\geq 4$ where the group $G\subset \SAut (X)$ is generated
by a saturated set $\cN\subset \LND(X)$ such that $0\ne \delta \in \cN$ and the next two  
conditions are true:

\noindent ($\#$) the variety $\rho^{-1} (Q\setminus Q_0)$ is  of codimension at least $2$ in $X$
(in particular, the zero locus of $\delta$ is at least of codimension 2);

\noindent ($\# \#$) $X$ is affine and parametrically $G$-suitable with respect to $\delta$.
\econv

\brem\label{nor.r1}  The condition  ($\#$) cannot be acheived for some flexible Gizatullin surfaces but in
higher dimensions the authors do no know  examples of smooth affine flexible varieties when  ($\#$) does not hold for an appropriate $\delta$.
For instance, if $X$ is a linear algebraic group without nontrivial characters, then every nonzero nilpotent element of its Lie algebra
yields a locally nilpotent vector field on $X$ for which condition ($\#$) is true.
Another example is a complex sphere given by $z_1^2+\ldots +z_m^2=1$ where $\bar z= (z_1, \ldots , z_m)$ is a coordinate system on $\A^m$. 
Then $\sigma =z_2\frac{\p}{\p z_1} - z_1\frac{\p}{\p z_2}$ is a locally nilpotent vector field on this sphere whose zero locus
has codimension 2.  The kernel of $\sigma$ is generated by $z_3, \ldots , z_m$ and, therefore, the quotient space is
smooth and ($\#$) holds.
\erem

\section{Normal bundles in flexible varieties}


Before describing $N_XC$ from Setting \ref{nor.s1} we need the following.

\bprop\label{nor.p1} Let   $X$ be $G$-flexible
where the group $G\subset \SAut (X)$ is generated
by a saturated set $\cN\subset \LND(X)$ such that $\delta \in \cN$.
Let  $K$ be a finite subset of $C$.
Then there exists a perfect $G$-family $\cA$ such that for general $\alpha_1, \ldots, \alpha_{n-1} \in \cA$
there is an open neighborhood $C^*$ of $K$ in $C$ for which the following is true.

{\rm (i)} For every $i=1, \ldots, n-1$ the restriction of $\rho_{\alpha_i}$ yields an isomorphism between $C^*$ and $C^*_i:=\rho_{\alpha_i} (C^*)$
where $C_i^*$ is contained in $Q_0$.

{\rm (ii)} The normal bundle $N_XC^*$ is trivial and $\pr (\delta_{\alpha_1}), \ldots, \pr (\delta_{\alpha_{n-1}})$ generate $N_XC^*$ as a module over $\bk [C^*]$.

{\rm (iii)}  The morphism $\cA \to \prod_{x \in K} T_x X, \, \alpha \mapsto ( \delta_\alpha(x) )_{x \in K}$ is dominant.
\eprop

\bproof  Statement (i) follows from Theorem \ref{pre.t4} and Remark \ref{pre.r3} (2) applied to $X_{\rm reg}$ and the fact that $\rho_\alpha =\rho\circ \alpha^{-1}$.
Let $T'X_{\rm reg}$ be the complement to the zero section in  the tangent bundle $TX_{\rm reg}$. Recall that $\cA$ acts transitively on $T'X_{\rm reg}$ (see, Example \ref{pre.exa1}).
Hence, letting $Y=T'X_{\rm reg}$  in Remark \ref{pre.r1} we get (iii).
In particular, we can suppose that the vectors
$\delta_{\alpha_1}(x), \ldots, \delta_{\alpha_{n-1}}(x)$ are general in $T_xX$ and, thus, linearly independent modulo $T_xC$.
Taking, if necessary, a  smaller neighborhood $C^*$ of $K$ such that for every $x \in C^*$ the similar fact hold we get (ii). Hence, we are done.
\eproof

It turns out that under Convention \ref{nor.conv1} we can get a much stronger version of Proposition \ref{nor.p1}.




\bprop\label{nor.p2} Let Convention \ref{nor.conv1} hold. 
Then one can choose a perfect $G$-family $\cA$ so
that for a general $\alpha \in \cA$ 

{\rm (i)} the morphism $\rho_\alpha|_C : C \to Q$ is proper;

{\rm (ii)}  its restriction yields an isomorphism
between $C$ and $\rho_\alpha (C)$ and 

{\rm (iii)} $\rho_\alpha (C) \subset Q_0$
(in particular,  $\delta_\alpha$ does not vanish on $C$).
\eprop

\bproof  
Theorem \ref{pre.t1} implies that for a general $\alpha$ in any perfect family  $\alpha^{-1} (C)$ does not meet $\rho^{-1} (Q\setminus Q_0)$,
since the latter is of codimension at least 2.  Since $\rho_\alpha=\rho \circ \alpha^{-1}$, this yields (iii).
Let $\cA_1$ now be a perfect $G$-family as in Theorem \ref{pre.t2}(ii), i.e. for general $\alpha\in \cA_1$ the morphism
 $\rho_\alpha|_C : C \to Q$ is an injective immersion.
Note that if $C$ is once-punctured, then $\rho_\alpha|_C : C \to Q$ is automatically proper and, in particular,
$\rho_\alpha (C)$ is closed in $Q$. 
Hence, in this case $\rho_\alpha|_C$ yields an isomorphism between $C$ and $\rho_\alpha (C)$.
If $C$ is not once-punctured consider a family $\cB\subset G$ of automorphisms of $X$ as in Definition \ref{prop.d2}.
Recall that by Remark \ref{pre.r2} the replacement of $\cA_1$ by $\cA=\cA_1 \times \cB$ leaves the conclusions of
Theorem \ref{pre.t2} valid.  Consider the morphism $\cA_1\times \cB \times X \to \cA_1\times \cB \times Q, \,
(\alpha, \beta ,x) \mapsto (\alpha, \beta, \rho_{\alpha\beta}(x))$ and extend it to a proper morphism
$\cX \to  \cA_1\times \cB \times Q$ over $ \cA_1\times \cB$. Let $\cD = \cX \setminus (\cA_1\times \cB \times X)$ and
$\cC$ be the closure of $\cA_1\times \cB \times C$ in $\cX$. 

Assume that  $\rho_{\alpha\beta}|_C : C \to Q$  is not proper  for a general $(\alpha, \beta)$ in $\cA_1\times \cB$.
Recall that 
the morphism $\rho_{\alpha\beta}|_C : C \to Q$  is not proper iff the closure $C_{\alpha, \beta}$ of the curve $\alpha\times \beta \times C$ in $\cX$ meets $\cD$.
Hence,  there exists  
a  dense subset $U$ of $\cA_1\times \cB$ such that for every $(\alpha, \beta ) \in U$ the curve $C_{\alpha, \beta}$ meets $\cC\cap \cD$.
This implies that there exists an irreducible component of $\cC \cap \cD$ whose image in $\cA_1 \times \cB$ under the natural morphism
$\cD \to \cA_1 \times \cB$ is dense.
Therefore, we can suppose that $U$ is open.
For a general $\alpha \in \cA_1$  and the natural projection
$\theta : \cA_1\times \cB \to \cA_1$ the intersection
of $\theta^{-1} (\alpha ) \cap U$ is open and nonempty in $\theta^{-1}(\alpha) \simeq \cB$.
That is, the morphism $\rho_{\alpha\beta}|_C : C \to Q$  is not proper when $\beta \in \theta^{-1} (\alpha ) \cap U$.
On the other hand,  
this morphism factors through the isomorphism $C \to \alpha^{-1} (C)$
and the map $\rho_\beta : \alpha^{-1} (C) \to Q$. By Theorem \ref{prop.t1} and the parametric suitability assumption
the latter map is proper for a general $\beta \in \cB$. This contradiction shows that our assumption is wrong which yields the desired conclusion.
\eproof

\blem\label{nor.l1} Let Convention \ref{nor.conv1}  hold and
$\cA$ be a perfect $G$-family as in Proposition \ref{nor.p2}.
Let $Z_\alpha$ be the zero locus of $\delta_\alpha$ and  $L_\alpha$ be the line bundle induced by $\delta_\alpha$ on $X\setminus Z_\alpha$.
Suppose that  $K$ is a finite subset of $C$ and $X_k$ is an open subset of $X$ containing $C$.
Let  $B_k\subset TX_k$ be a vector subbundle of rank $k \leq n-2$. For general   $\alpha$ in $\cA$ one has the following.

{\rm (i)}  There exists a  closed subvariety  $Y_k\subset X$ of dimension at most $k$ such that 
$B_k$ and $L_\alpha$ generate a vector bundle  $B_{k+1}$ of rank $k+1$ on  $X_{k+1}=X_k\setminus  (Y_k\cup Z_\alpha)$
and $C\subset X_{k+1}$.

{\rm (ii)} If  $k\leq n-3$ and $\pr (B_k)$ yields a vector subbundle of $N_X C$ of rank $k$,  then $\pr (B_{k+1})$
 yields a vector subbundle of $N_XC$ of rank $k+1$.
 
 {\rm (iii)}  If  $\pr (B_{n-2})$ is a vector subbundle of $N_X C$ of rank $n-2$, then  there exists a finite set
 $H_{\alpha} \subset  C\setminus K$ such that $\pr (B_{n-1})$
 yields $N_X( C\setminus H_{\alpha})$.  Furthermore,  $\bigcap_{\alpha \in \cA}H_{\alpha}= \emptyset$.

\elem

\bproof 
Recall that by $G$-flexibility and \cite[Theorem 4.14]{AFKKZ} $G$ acts transitively on $T'X$ where $T'X$
is the complement to the zero section in $TX$.
For every subset $S$ of $TX$ let $S'=S\cap T'X$.  Note that $\dim B_k'=n+k$, $\dim B_k'|_C=k+1$ and $\dim L_\alpha' =n+1$.
By Theorem \ref{pre.t1} for general $\alpha \in \cA$ the variety $L_\alpha'$  meets  $ B_k'|_C$ along a subvariety $R$ of dimension
$(k+1)+(n+1)-2n=k+2-n\leq 0$. However, this dimension cannot be zero, since
for each $v\in R\subset TX$ every vector proportional to $v$ is also in $R$, i.e., $R$ is empty.
Again by Theorem \ref{pre.t1}  $L_\alpha'$ meets $B_k'$ transversely along a subvariety $\tY$
of dimension at most $(n+k)+(n+1)-2n=k+1$. Hence,  $Y_k=\pr (\tY)$ does not meet $C$ and has dimension at most $\dim \tY -1=k$,
while $B_k$ and $L_\alpha$ generate a rank $k+1$ vector subbundle $B_{k+1}$ of the tangent bundle of   $X_{k+1}$.
Note that $C\subset X_{k+1}$ since $C\cap Z_\alpha =\emptyset $ by Proposition \ref{nor.p2} (iii). Thus, we get (i).

The claim in (ii) that $\pr (B_k)$ is a vector subbundle of $N_X C$ of rank $k$  is equivalent to
the fact  that the constructible set $\bigcup_{x\in C}(B_k(x)+T_xC)$  is a bundle $\tB_k$ of rank $k+1$ over $C$.  
We note that since  $\dim \tB_k'+\dim L_\alpha' -\dim T'X=(k+2)+(n+1)-2n=k+3-n\leq 0$, by Theorem \ref{pre.t1} for general $\alpha \in \cA$ 
the variety $L_\alpha'$ meets $\tB_k'$ transversely along a variety $\tH$ of dimension at most zero. However, arguing as before we 
see that $\tH=\emptyset$ which implies (ii) (to see that each $x \in K$ is still in $X_{k+1}$ one needs to choose as before
$\delta_\alpha (x)$ outside the fiber $B_k(x)+T_xC$).

In (iii) the same argument shows that $\tH$ may have dimension 1.  Therefore, $H_\alpha :=\pr (\tH)$ is a finite set
and $\pr (B_{n-1})$ yields $N_X( C\setminus H_{\alpha})$.
The transversality argument  implies also that
for general $\alpha \in \cA$ the variety $L_\alpha'$ does not meet the fiber of $\tB_{n-s-1}'$ over a fixed point $c \in C$.
This implies that $\bigcap_{\alpha \in \cA}H_\alpha= \emptyset$ and, consequently,
$K\cap H_\alpha=\emptyset$ for general $\alpha$ which concludes the proof.
\eproof

\bprop\label{nor.p3} Let Convention \ref{nor.conv1} hold
and  $K$ be a finite subset of  $C$.
Then  there exists a perfect $G$-family $\cA$ 
such that  for $\bar \alpha =(\alpha_1, \ldots, \alpha_{n-1})$ in a dense open subset $U$ of $\cA^{n-1}$
the following is true.

 {\rm (1}) For every $i=1, \ldots ,n-1$ the map $\rho_{\alpha_i} : C \to Q$ is proper and generates an isomorphism
 between $C$ and $\rho_{\alpha_i} (C)\subset Q_0$.

{\rm (2)} The span of $\delta_{\alpha_1}, \ldots ,\delta_{\alpha_{n-2}}$ over $\bk [C]$ is a subbundle $B_{n-2}\subset TX|_C$  such that
$\pr (B_{n-2})$ is a subbundle of $N_XC$ whose rank is $n-2$. 

{\rm (3)} Let $H_{\alpha_n}$ be as in Lemma \ref{nor.l1} and $C^*=C\setminus H_{\alpha_{n-1}}$.
Then $K\subset C^*$ and the span of $\pr (B_{n-2})$ and
 $\pr(\delta_{\alpha_{n-1}})$ over $\bk [C^*]$ is the normal bundle $N_X C^*$.

{\rm (4)}  The morphism $\cA \to \prod_{x \in K} T_x X, \, \alpha \mapsto ( \delta_\alpha(x) )_{x \in K}$ is dominant.

{\rm (5)} $N_XC$ is naturally isomorphic to $\pr (B_{n-2})\oplus L$ where $L$ is a line subbundle of $N_XC$. 
 \eprop

\bproof  The existence of a perfect $G$-family $\cA$ satisfying (1)  is provided by Proposition \ref{nor.p2}.
Let $B_0$ be the zero subbundle of $TX$. Then applying Lemma \ref{nor.l1} to $B_0$ and consequent vector bundles $B_k$
we get (2) and (3). The same argument as in Proposition \ref{nor.p1} (iii) yields (4).

Now let  $v_i=\pr (\delta_{\alpha_i})$.  
Then each $v_i$ is a non-vanishing section of $N_XC$ and by (3) $v_{n-1}$ meets $\pr (B_{n-2})$
over $H_{\alpha_{n-1}}=\{ x_1, \ldots , x_m\}\subset C$.  Consider a neighborhood $U_i$ of $x_i$ in $C$ such that $U_i\cap H_{\alpha_{n-1}}=x_i$ and
the normal bundle $N_XU_i$ is trivial, i.e., $v_1, \ldots, v_{n-2}$ and some section $u_i$ is a basis of $N_XU_i$.
Then $v_{n-1}=g_iu_i+f_{1,i}v_1+ \ldots + f_{n-2,i} v_{n-2}$ where $g_i,f_{1,i}, \ldots, f_{n-2,i}$ are regular functions
on $U_i$.  
Consider regular functions
$f_1, \ldots,  f_{n-2}$ on $C$ such that $f_j-f_{j,i}$ is divisible by $g_i$ for every $i$ and $j$.
Consider also the projectivization $D$ of $N_XC$ which is a $\PP^{n-2}$-bundle over $C$.
It contains a $\PP^{n-3}$-bundle $D_0$ corresponding to $\pr (B_{n-2})$. 
Note that $v_{n-1}- (f_{1}v_1+ \ldots + f_{n-2} v_{n-2})$ defines a section
of $D$ over $C\setminus H_{\alpha_{n-1}}$
which extends to a section over $C$ that does not meet $D_0$. This section yields the desired line bundle $L$ which concludes the proof.
\eproof

\brem\label{tech.r1}    Recall that by \cite{Se1} every vector bundle $W$ of rank $k$ over $C$ can be presented as $W_0\oplus W_1$ where
$W_0$ is a trivial bundle of rank $k-1$ and $W_1$ is a line bundle.  
Proposition \ref{nor.p3}  (5) below tells additionally that for $W=N_XC$ the bundle $W_0$ has a basis of sections generated
by locally nilpotent vector fields.


\erem

\bcor\label{nor.c1} Let the notations and conclusions of Proposition \ref{nor.p3} hold and $v_i=\pr (\delta_{\alpha_i}), \, i=1, \ldots , n-1$. 
Suppose that $N_CX$ and, therefore, $L$ are trivial bundles and $v_0$ is a nonvanishing section of $L$.
Then $v_{n-1} = q v_0+\sum_{i=1}^{n-2} a_i v_i$ where $a_1, \ldots , a_{n-2}$ and $q$ are regular functions on $C$
 having  simple zeros on $C$.
\ecor

\bproof  Consider the natural projection $N_XC \to L$. Let $\tilde v_0=qv_0$ be the image of $v_{n-1}$.
Then $v_{n-1}-\tilde v_0\in \pr (B_{n-2})$ which yields an equality $v_{n-1} =\tilde v_0+\sum_{i=1}^{n-2} a_i v_i$. 
Since $\pr (B_{n-2})$ is a vector subbundle of $N_X C$ of rank $n-2$,
the constructible set $\bigcup_{x\in C}(B_{n-2}(x)+T_xC)$  is a bundle $\tB_{n-2}$ of rank $n-1$ over $C$.
Let $L_{\alpha_{n-1}}$  be the line bundle induced by $\delta_{\alpha_{n-1}}$ on $X\setminus Z_{\alpha_{n-1}}$
as in Lemma \ref{nor.l1}.
As we mentioned in the proof of  Lemma \ref{nor.l1} $L_{\alpha_{n-1}}$ meets $\tB_{n-2}$ transversely 
for a general $\alpha_{n-1} \in \cA$ by virtue of Theorem \ref{pre.t1}
which implies the claim about the zeros of $q$ (which are the same as the zeros of $\tilde v_0$).
The similar argument works for the zeros of $a_i$. Thus, we get the  desired conclusion.
\eproof

\section{$\SL_n (A)$ over a commutative ring  $A$}

In this section we collect some facts about matrices over commutative rings and prove Proposition \ref{bass.p1} which is an essential
step for our main result.

\bnota\label{bass.n1} Let  $A$ be a commutative ring,
$n \geq 3$ and $\SL_n (A)$ be the special linear group with entries from $A$,
whereas ${\rm Mat}_n (A)$ is the ring of $(n\times n)$-matrices with entries from $A$. 
By $I\in \SL_n (A)$ we denote the identity matrix and by $e_{ij}, \, 1\leq i, j\leq n$ the matrix in ${\rm Mat}_n (A)$ 
with all zero entries but the entry in position $(i,j)$ which is 1.
By $\EE_n (A)$ we denote the subgroup of $\SL_n(A)$ generated by all elementary matrices of the form $I +a e_{ij},\, 1\leq i\ne j\leq n$ where $a \in A$.
Let $\fq$ be any ideal in $A$.  
By  $\SL_n (A, \fq)$ we denote
the kernel of the natural group homomorphism  $\SL_n(A) \to \SL_n (A/\fq)$, 
by $\FF_n(A,\fq)$ we denote the subgroup of $\EE_n(A)$
generated by all elementary transformations of the form $I +a e_{ij}$ where $a \in \fq$
and by $\EE_n (A, \fq)$ the normal closure of $\FF_n(A, \fq)$ in $\EE_n (A)$. 
\enota

We need the following result of Bass, Milnor and Serre \cite[Theorem 4.1]{BMS} (see also \cite[Theorem 4]{Ba}).

\bthm\label{bass.t1} 
Let Notation \ref{bass.n1} hold, $A$ be a Dedekind  domain and $n\geq 3$. Then
 \be\label{bass.eq1}{\rm E}_n (A,\fq)=[{\SL}_n (A), {\SL}_n (A,\fq)].\ee 
\ethm

Another important fact unfortunately omitted in \cite{Ud} is the following ``Tits' theorem" in \cite[Theorem 3]{Ni}.

\bthm\label{bass.t2} Let Notation \ref{bass.n1} hold. Then $\EE_n(A,\fq^2)$ is contained in $\FF_n(A,\fq)$.

\ethm

Let us reproduce two facts from \cite[Lemma 2.2.1 and Theorem 2.2.2]{Ud}.

\blem\label{bass.l1} {\rm (\cite[Theorem 2.2.2]{Ud})} Let Notations \ref{bass.n1} hold, $\fq=qA, \, q \in A$ be a proper principal ideal in $A$
and $D\in {\rm Mat}_n(A)$ be a matrix which is invertible in ${\rm Mat}(A[\frac{1}{q}])$.

{\rm (i)}  For some $m >0$ one has $ D^{-1} {\SL}_n(A, \fq)D\supset {\SL}_n(A, \fq^{m})$.

{\rm (ii}) Let $A$ be a Dedekind ring. Then the group $D^{-1} E_n(A, \fq)D$ contains $E_n(A, \fq^l)$ for some $l>0$.

\elem

\bproof Since $D^{-1} =\frac{1}{{\rm det} D}{{\rm adj} (D)}$ (where  ${\rm adj} (D)$ is the classical adjoint),  one has (i)
with $m$ such that $\frac{q^{m-1}}{{\rm det} D}$ belongs to $A$.

By (i) and  Formula \eqref{bass.eq1} we have
$$D^{-1} {\rm E}_n(A, \fq)D=D^{-1} [{\SL}_n(A),{\SL}_n(A, \fq)]D=[D^{-1} {\SL}_n(A)D, D^{-1}{\SL}_n(A, \fq)D] \supset $$
$$[ {\SL}_n(A, \fq^m), {\SL}_n(A, \fq^m)] \supset [ {\rm E}_n(A, \fq^m), {\rm E}_n(A, \fq^m)]\supset {\rm E}_n(A, \fq^{2m}) $$
where the last inclusion follows from \cite[Section 5]{BMS}.
\eproof 

Now we have the desired fact.

\bprop\label{bass.p1} Let Notation \ref{bass.n1} hold, $A$ be a Dedekind ring and let $q$ and $q' \in A$ generate comaximal principal ideals $\fq$ and $\fq' \subset B$.
Let a matrix $D'$  plays the same role for $q'$ as $D$ for $q$  in Lemma \ref{bass.l1}.
Suppose that $H$ is a subgroup of $\SL_n (A)$ that contains all elements from
$\SL_n(A)\cap (D^{-1} \FF_n(A, \fq)D)$ and $\SL_n(A)\cap ((D')^{-1} \FF_n(A, \fq')D')$.
Then  $H$ contains $\EE_n (A)$.
\eprop

\bproof  By Theorem \ref{bass.t2}, $H$ contains the subgroups
$\SL_n(A)\cap (D^{-1} \EE_n(A, \fq^2)D)$ and $\SL_n(A)\cap ((D')^{-1} \EE_n(A, (\fq')^2)D')$.
Hence, by Lemma \ref{bass.l1} $H$ contains $\EE_n(A, \fq^{2l})$ and $\EE_n(A, (\fq')^{2l})$.
Now the desired conclusion follows from the fact that  $(\fq')^{2l}+\fq^{2l} =A$.
\eproof

The rest of this section will be needed only for the holomorphic case in Section 7.

\blem\label{bass.l2} 
Let Notation \ref{bass.n1} hold, $\sigma \in {\rm Mat}_n (A)$ be a matrix of the form ${q}(I +\sum_{i=2}^n \frac{a_i}{q}e_{1i})$ 
where $ q \in A\setminus \{ 0\}$ and each $a_i\in A$.  Suppose that $\fp$ is an ideal in $A$.
Let $\cE$ be a subgroup of $\SL_n (A)$ that contains 

{\rm (i)} all elementary matrices of the form $\tau =I+ce_{ij}$ where  $c\in A$, 
$2\leq j\ne i \leq n$ and

{\rm (ii)} all  matrices $\sigma \tau \sigma^{-1}$ where $\tau =I+bq^2e_{i1}$, $b\in \fp$  and $1\leq i \leq n$.

Then  for every $2 \leq k\ne l \leq n$ and every $b \in \fp$ the group $\cE$ contains the following matrices
$$I+ a_lbq^2e_{k1}  -a_ka_lbqe_{kk} -a_l^2bqe_{kl} - a_kbq^2e_{l1} + a_ka_lbqe_{ll}  + a_k^2bqe_{lk}.$$
\elem

\bproof   Note that $\sigma^{-1} = {q}^{-1}(I -\sum_{i=2}^n \frac{a_i}{q}e_{1i})\in {\rm Mat}_n ({\rm Frac} (A))$.
Let $\alpha =I+ba_lq^2e_{k1}-ba_kq^2e_{l1}$. Then
$\sigma \alpha \sigma^{-1}= I + \sigma(ba_lq^2e_{k1}-ba_kq^2e_{l1})\sigma^{-1}=
I+ \sigma(ba_lq^2e_{k1})\sigma^{-1} -   \sigma(ba_kq^2e_{l1})\sigma^{-1}$. Consider
\[ \begin{array}{c} \sigma(bq^2e_{k1})\sigma^{-1} = bq^2e_{k1} +
(\sum_{i=2}^n \frac{a_i}{q}e_{1i}) bq^2e_{k1}- bq^2e_{k1}(\sum_{i=2}^n \frac{a_i}{q}e_{1i})\\
-(\sum_{i=2}^n \frac{a_i}{q}e_{1i}) bq^2e_{k1}(\sum_{i=2}^n \frac{a_i}{q}e_{1i})=\\
bq^2e_{k1} +a_kbqe_{11} -a_kbqe_{kk} -a_k^2be_{1k}
- \sum_{i=2, i\ne k}^n a_ibqe_{ki} -\sum_{i=2, i\ne k}^n a_ka_ib e_{1i}=\\
bq^2e_{k1} +a_kbqe_{11} -a_kbqe_{kk} -a_k^2be_{1k} - a_lbqe_{kl}- a_ka_lb e_{1l}\\
-\sum_{i=2, i\ne k,l}^n a_ibqe_{ki} -\sum_{i=2, i\ne k,l}^n a_ka_ib e_{1i}.\end{array}\]
Replacing $b$ by $a_lb$ we get  $\sigma (ba_lq^2e_{k1})\sigma^{-1}$ and then, exchanging the letters $k$ and $l$ 
we get $ \sigma (ba_kq^2e_{l1})\sigma^{-1}$. Hence, one has
\[ \begin{array}{c} \sigma \alpha \sigma^{-1}=
I+a_lbq^2e_{k1} +a_ka_lbqe_{11} -a_ka_lbqe_{kk} -a_k^2a_lbe_{1k} - a_l^2bqe_{kl}- a_ka_l^2b e_{1l}\\
-( a_kbq^2e_{l1} +a_ka_lbqe_{11} -a_ka_lbqe_{ll} -a_ka_l^2be_{1l} - a_k^2bqe_{lk}- a_k^2a_lb e_{1k}) +\beta\\
= I+ a_lbq^2e_{k1}  - a_ka_lbqe_{kk} - a_l^2bqe_{kl} - a_kbq^2e_{l1} + a_ka_lbqe_{ll}  + a_k^2bqe_{lk} +\beta                        \end{array}\]
where $\beta$ is a sum of terms of the form $\sum_{i=2, i\ne k,l}^n (c_ie_{ki} +d_ie_{li})$ where $c_i,d_i \in A$.
One can rewrite $ \sigma \alpha \sigma^{-1}$ now as
$$(I +\beta)(I+ a_lbq^2e_{k1}  -a_ka_lbqe_{kk} - a_l^2bqe_{kl} - a_kbq^2e_{l1} +  a_ka_lbqe_{ll}  + a_k^2bqe_{lk}).$$
Since $I+\beta \in \cE$, we have the desired conclusion.
\eproof

\bexa\label{bass.exa1}  Consdier the case of $(k,l) = (2,3)$ in the proof of Lemma \ref{bass.l2}. 
Let  $ \sigma \alpha \sigma^{-1}$ and $\beta$ be as in the proof of Lemma \ref{bass.l2} and $ \theta=(I+\beta)^{-1} \sigma \alpha \sigma^{-1}$. Then
\[ \theta =\left[ {\begin{array}{cccc}
1 & 0 & 0 & \bar 0\\
 a_3bq^2&  1-a_2a_3bq & -a_3^2bq & \bar 0  \\      
- a_2bq^2 & a_2^2bq & 1+ a_2a_3bq  &  \bar 0    \\
\bar  0 & \bar 0   & \bar 0 &  I_{n-3}  \end{array} } \right] \]
where $I_{n-3}\in \SL_{n-3}(A)$ is the identity matrix.
Removing the first row and the first column we get the matrix
\[   \gamma =\left[ {\begin{array}{ccc}
 1-a_2a_3bq & -a_3^2bq & \bar 0  \\      
 a_2^2bq & 1+ a_2a_3bq  &  \bar 0    \\
 \bar 0   & \bar 0 &  I_{n-3}  \end{array} } \right] \]
which has determinant 1.
\eexa

\bprop\label{bass.p2} Let the assumptions of Lemma \ref{bass.l2} hold and $\SL_{n-1}(A)=\EE_{n-1} (A)$.
Then  for every $b \in \fp$ the matrix $I + a_lbq^2e_{k1}  - a_kbq^2e_{l1}$ is contained in $\cE$.
\eprop
\bproof  We show only  that  $I + a_3bq^2e_{21} -a_2bq^2e_{31}\in \cE$ since the general case of $(k,l)\ne (2,3)$ is similar. 
Hence, we use the notations $\theta$ and $\gamma$ as in Example \ref{bass.exa1}.
Consider a block  matrix of the form
\[  \bar \lambda =\left[ {\begin{array}{cc}
 1 & \bar 0  \\      
  \bar 0   & \lambda  \end{array} }\right] \in {\rm SL}_n(A).\]
Then 
\[  \theta \bar \lambda =\left[ {\begin{array}{cccc}
1 &  \bar 0\\
 \bar c &\gamma \lambda  \end{array} }\right]\]
where $\bar c$ is an (n-1)-column with the first two entries $a_3bq^2$ and  $ -a_2bq^2$,
whereas the rest of the entries are zeros. Thus, taking $\lambda =\gamma^{-1}\in \SL_{n-1}(A)=\EE_{n-1}(A)$
we see that $I+ a_3bq^2e_{21}  - a_2bq^2e_{31}\in \cE$ which concludes the proof.
\eproof

\section{Main Theorem}

In this section we present a strengthened version of the main result in \cite{Ud}. Let us start with terminology.

\bdefi\label{main.d1} (1) Let $C_1$ and $C_2$ be closed subvarieties of a smooth quasi-affine variety $X$ with defining ideals $I_1$
and $I_2$ in $\bk [X]$.  Let $G$ be a subgroup of $\Aut (X)$, i.e., $G$ acts on $\bk [X]$. Suppose that for some $\varphi \in G$ one has $\varphi (I_1)=I_2$.
Then $\varphi (I_1^k)=I_2^k$ for every $k \in \N$ and, therefore, $\varphi$ induces an isomorphism  $\varphi_k : \frac{\bk [X]}{I_1^k} \to \frac{\bk [X]}{I_2^k}$.
We say that such a $\varphi_k$ is a {\em $G$-induced isomorphism}.

(2)   Let $\theta : \frac{\bk [X]}{I_1^k} \to \frac{\bk [X]}{I_2^k}$ be a ring homomorphism. Note that one has $\theta ( I_1) \subset {I_2}$
since $\theta$ sends nilpotent elements to nilpotent ones. Hence, $\theta$ induces ring  homomorphisms $\theta' : \frac{\bk [X]}{I_1} \to \frac{\bk [X]}{I_2}$
and   $\theta'' : \frac{I_1}{I_1^2} \to \frac{I_2}{I_2^2}$. 
Given an isomorphism $\lambda : \frac{I_1}{I_1^2} \to \frac{I_2}{I_2^2}$ we say that it is {\em $G$-induced} if for some $G$-induced isomorphism
$\theta : \frac{\bk [X]}{I_1^k} \to \frac{\bk [X]}{I_2^k}$ one has $\lambda =\theta''$.

(3) If in the previous definition  $I_1=I_2$ and  $\theta'$ is the identity map, then we say $\theta$ (resp. $\lambda$) is an {\em $A_1$-automorphism} of 
$\frac{\bk[X]}{I_1^k}$ (resp. $\frac{I_1}{I_1^k}$)
where $A_j=\frac{\bk [X]}{I_j}$ is the ring $\bk [C_j]$ of regular functions on $C_j$.  That is, $\theta$ can be viewed as
an automorphism of an $k$th infinitesimal neighborhood  of $C_1$ over $C_1$ itself.
\edefi

In the case of a flexible variety the Jacobian of an automorphism must be a constant function because of the following. 

\blem\label{main.l1} Let $X$ be a smooth quasi-affine algebraic variety equipped with a volume form $\omega$ such that
every two points of $X$ can be joined by a chain of polynomial curves (which is the case when $X$ is flexible). Let $\alpha :X \to X$ be an automorphism.
Then the Jacobian of $\alpha$ is constant. Furthermore, if $\alpha\in \SAut (X)$ then its Jacobian is 1.
\elem

\bproof The action of $\alpha$
transforms $\omega$ into a volume form $\omega_\alpha$. The Jacobian $f_\alpha=\frac{\omega_\alpha}{\omega}$ of this automorphism is an invertible function on $X$
and $f_e= 1$ where $e$ is the identity automorphism.  Since  the restriction of $f_\alpha$ to every polynomial curve is constant,
we see that $f_\alpha$ is constant as $\bk$ is an algebraically closed field.
Thus, for $\alpha \in \SAut (X)$ the Jacobian $f_\alpha$ can be viewed
as a value of an invertible function  $f$ on $\SAut (X)$. By the definition any two elements $\alpha_1,\alpha_2 \in \SAut (X)$ can be joined by
a chain of  polynomial curves (corresponding to one-parameter unipotent groups). Hence, $f$ is constant and, consequently, $f=1$ which is the desired conclusion.
\eproof

To see how the above condition on Jacobians can be reformulated for automorphisms of infinitesimal neighborhoods of closed subvarieties of $X$
we need to remind the following well-known fact.

\blem\label{main.l2} Let $Y$ be a normal affine algebraic variety of dimension $m$ and
$Z$ be a smooth closed curve in $Y$ with  defining ideal $J$ in $\bk [Y]$.  Suppose that $Z$ is contained in $Y_{\rm reg}$. Then any section of
the conormal bundle of $Z$  in $Y$ can be presented by a regular function on $Y$. Furthermore, suppose that
the  conormal bundle is trivial.  
Then
there exists an open neighborhood $U$ of $Z$ in $Y$ and regular functions $u_1, \ldots , u_{m-1} \in \bk [Y]$ such that $Z$ is a strict
complete intersection in $U$ given by $u_1|_U= \ldots = u_{m-1}|_U=0$.
\elem
\bproof Let $\cJ$ be the ideal sheaf on $Y$ induced by $J$. Consider the short exact sequence
$ 0\to \cJ^2 \to \cJ\to \frac{\cJ}{\cJ^2}\to 0$ and the long exact sequence induced by it. Since $\cJ^2$ is a coherent sheaf  on the affine variety $Y$,
we have $H^1(Y, \cJ^2)=0$ and the long exact sequence implies that every section of $\frac{\cJ}{\cJ^2}$ is induced by a section of $\cJ$, i.e.,
by a regular function on $Y$ which belongs to $J$. 
Note also that as in \cite[Chapter II,  Theorem 8.17]{Har} $ \frac{\cJ}{\cJ^2}$ can be viewed as 
a sheaf $\cC$ of $\cO_Z$-modules where $\cO_Z$ is the structure sheaf of $Z$.
Any section of $\cC$ over $Z$ can be extended by 0 to $Y\setminus Z$ so that one gets a global section of $ \frac{\cJ}{\cJ^2}$
which, as we showed, is induced by an element of $J$.
Furthermore, since $Z$ is smooth   $\cC$
is  a locally free sheaf of rank $m-1$ 
which induces the conormal bundle of $Z$ in $Y$.
Since any section of this bundle can be viewed as a section of $\cC$, we have the first statement.
The smoothness of $Z$ implies that $Z$
 is a local complete intersection in $Y$ (\cite[Chapter II, Theorem 8.17]{Har}). That is, for every point $z \in Z$ and any collection $u_1^z, \ldots, u_{m-1}^z\in \cJ$
which induces a local basis $\tilde u_1^z, \ldots, \tilde u_{m-1}^z$ of $\cC$ at $z$ one can find a neighborhood $U_z$ of $z$ in $Y$
such that $U_z \cap Z$ is  a strict complete intersection in $U_z$ given by  the system $u_1^z= \ldots = u_{m-1}^z=0$. 
In the case of a trivial $\cC$ for every $i=1, \ldots, m-1$ such local sections  $\tilde u_i^z$ can be chosen as the restrictions of a global section $\tilde u_i$ of $\cC$ to $Z\cap U_z$.
As we saw $\tilde u_i$ is induced by a function $u_i \in J$. 
There is  also a neighborhood $U$ of $Z$ in $Y$, for which $U\cap Z$ is given by the system  $u_1= \ldots = u_{m-1}=0$,
since for every $z \in Z$ there is a neighborhood $W$ of $z$ in $Y$ in which $W\cap Z$ is given by this system.
This concludes the proof.
\eproof

\bdefi\label{main.d2}
Let  $X$ be a normal quasi-affine algebraic variety, i.e., $X$ is contained as a dense subvariety in a normal affine variety $Y$.
Let $C_j$ be a smooth curve in $X_{\rm reg}$ with a trivial normal bundle such that  $C_j$ is closed in $Y$.
 By Lemma \ref{main.l2}
$C_j$ is a strict complete intersection in some open $U_j\subset X$ given
by $u_{1j}= \ldots = u_{n-1,j}=0$ where  $u_{1j}, \ldots , u_{n-1,j}$ are regular functions on $X$.
One can see now that the ring $\frac{\bk [X]}{I_j^k}$ is naturally isomorphic to $A_j\oplus \bigoplus_{l=1}^{k-1}\frac{I_j^l}{I_j^{{l+1}}}$ where each summand
$\frac{I_j^l}{I_j^{l+1}}$ can be treated as the space of homogeneous polynomials in $u_{1j}, \ldots , u_{n-1,j}$ of degree $l$ over $A_j$.
In particular,  for any homomorphism $\theta : \frac{\bk[X]}{I_1^k} \to \frac{\bk[X]}{I_2^k}$  of algebras one can consider the Jacobi
matrix $\left[\frac{\p \theta (u_{l1})}{\p u_{m2}}\right]_{l,m=1}^{n-1}$. In the case of $I_1=I_2$ and $u_{i,1}=u_{i,2}$
the determinant of this matrix is independent of the choice of coordinates $u_{1j}, \ldots , u_{n-1,j}$
and we say that  $\theta$ has Jacobian 1  if the determinant is equal to 1 modulo $I_1^{k-1}$.  
\edefi

Let  $\theta:  \frac{\bk[X]}{I_1^k} \to \frac{\bk [X]}{I_1^k}$ be  an $A_1$-automorphism induced by an automorphism $\alpha \in \Aut (X)$.
If $X$ possesses a volume form and the Jacobian of $\alpha$ is 1, then the Jacobian of $\theta$ is, of course, 1.
However, one can often make the same conclusion about the Jacobian of $\theta$ without the presence of a volume form.

\bexa\label{main.exa1}  Let $0\ne \delta \in \LND (X)$ and $\rho : X\to Q$ be a partial quotient morphism associated with $\delta$.
Suppose that $x\in C_1$ is  a general point of $X$, i.e., $\rho (x)$ belongs to $Q_0$ where $Q_0$ is as in Remark \ref{def.r1}.
Let $f\in \bk[Q]$ vanish on $\rho (C_1)$
and $\alpha\in \G_a$ be the automorphism of $X$ induced by $\exp( \rho^*(f)\delta)$. Then for every $k$ the automorphism
$\theta:  \frac{\bk[X]}{I_1^k} \to \frac{\bk [X]}{I_1^k}$ induced by $\alpha$ has Jacobian 1. 
Indeed, consider a function $h \in \bk[X]$ of degree 1 with respect to $\delta$, i.e., $h \notin \Ker \delta$ and $\delta(h) \in \Ker \delta$.
This is equivalent to the fact that the restriction of $h$ to the $\G_a$-orbit $\A^1$ of any general point of $X$ (and $x$, in particular) is a 
polynomial of degree 1. Consider the smooth germ $S$ of the hypersuface $h^{-1}(h(x)) \subset X$ at $x$. Then the $\G_a$-orbit of $S$ is naturally isomorphic
to $S\times \A^1$ and one can choose a coordinate $t$ on the second factor of $S\times \A^1$ so that
$\alpha|_{S\times \A^1}$ is given by $(s,t) \mapsto (s,t+f(s))$. 
Since $S$ is a smooth germ it has a volume form $\omega_0$ and, hence, $S\times \A^1$ has the volume form $\omega=\omega_0\times \dd t$.
It remains to note now that $\alpha$ preserve $\omega$ which yields the desired claim.  

\eexa

The main strengthened result  from \cite{Ud} is the following.

\bthm\label{main.t1} {\rm (\cite[Theorem 2.5.19]{Ud})} Let 
$X$ be a normal $G$-flexible quasi-affine variety of dimension $n\geq 4$ where $G\subset \SAut (X)$
is a group generated by a saturated set $\cN\subset \LND (X)$. Let  $0\ne \delta\in \cN$, 
$\rho : X \to Q$ be an associated partial quotient morphism  of $\delta$ and 
$C$ be a smooth  curve in $X$ such that $C\subset X_{\rm reg}$ and $C$ is closed in an affine variety containing $X$.  Suppose also that

{\rm ($\bullet $)} \hspace{0.5cm} $C$ is irreducible and either $C$ is once-punctured  or Convention \ref{nor.conv1} is valid
(in particular, $X$ is affine and parametrically suitable with respect to $\delta$).

\noindent Let  $I$ be the defining ideal of $C$ in the algebra $\bk [X]$ of regular functions on $X$ and $A=\frac{\bk [X]}{I} \simeq \bk [C]$.
Suppose that  the normal bundle $N_XC$ of $C$ in $X$ is trivial (i.e., it is a free $A$-module of dimension $n-1$).
Let $\theta : \frac{\bk[X]}{I^k} \to \frac{\bk [X]}{I^k}$ be  an $A$-automorphism with
Jacobian 1 (in particular, the automorphism $\lambda :  \frac{I}{I^2} \to \frac{I}{I^2}$ induced by $\theta$
belongs to $\SL_{n-1} (A)$).  Then one has the following

{\rm (i)} If $\lambda $ belongs to $\EE_{n-1} (A)$, then $\lambda$ is $G$-induced.

{\rm (ii)} If $\EE_{n-1} (A)=\SL_{n-1}(A)$ (which is true in the case when $A$ is a Euclidean domain \cite[Chapter I, Proposition 5.4]{Lam}, e.g., 
when $C$ is a smooth polynomial curve),
then $\theta$ is $G$-induced.

\ethm

Before the proof we need to consider several facts.  The first of which follows immediately from  Formula \eqref{pre.eq3} in Lemma \ref{pre.l1}.

\blem\label{main.l3} {\rm (cf. \cite[page 47]{Ud})}  Let $X, G, C$, $\delta$ and $\rho : X \to Q$ be as in Theorem \ref{main.t1}
and $\alpha_1, \ldots , \alpha_{n-1} \in G$. Let $\delta_\alpha =\alpha_* (\delta)$ and $\rho_\alpha =\rho\circ \alpha^{-1} : X \to Q$ for every $\alpha \in \Aut (X)$.
Suppose that $C_i:=\rho_{\alpha_i} (C)$ for some  $i\in \{ 1, \ldots,  n-1\}$. 
Suppose that $f_j$ and $h_j, \, j=1, \ldots, n-1$ are regular functions on $Q$ such that each $f_j$ vanishes on $C_i$.
Let $\Phi^{ij}_t$ be the flow of the locally nilpotent vector field $\rho_{\alpha_i}^*(h_jf_j) \delta_{\alpha_i}$ at time $t$
and let $\breve \sigma =\sigma|_C$ for every $\sigma \in \LND (X)$ .

{\rm (i)} 
Then for every $t$  the automorphism $\Phi^{ij}_t$ leaves each point of $C$ fixed (i.e., it preserves $TC$)  and it induces
an automorphism of the bundle $TX|_C$ (and, hence, of $N_XC$) such that $\breve \delta_{\alpha_i}\mapsto \breve \delta_{\alpha_i}$ and 
$\breve \delta_{\alpha_k}\mapsto \breve \delta_{\alpha_k}+t (\rho_{\alpha_i}^*(h_j)\delta_{\alpha_k}(\rho_{\alpha_i}^*(f_j)))|_C\breve \delta_{\alpha_i}$ for $1\leq k\ne i \leq  n-1$.

{\rm (ii)} For $j\ne l$ the locally nilpotent vector fields $\rho_{\alpha_i}^*(h_jf_j) \delta_{\alpha_i}$ and  $\rho_{\alpha_i}^*(h_lf_l) \delta_{\alpha_i}$ commute, 
i.e., the vector field $\sigma=(\sum_{j=1, j\ne i}^{n-1} \rho_{\alpha_i}^*(h_jf_j) )\delta_{\alpha_i}$
is locally nilpotent.  In particular, the flow of $\sigma$ at time $t$ induces an automorphism of $TX|_C$ such that 
$\breve \delta_{\alpha_i}\mapsto \breve\delta_{\alpha_i}$ and 
$$\breve \delta_{\alpha_k}\mapsto \breve \delta_{\alpha_k}+t (\sum_{j=1, j\ne i}^{n-1}  \rho_{\alpha_i}^*(h_j)\delta_{\alpha_k}(\rho_{\alpha_i}^*(f_j)))|_C 
\breve \delta_{\alpha_i}$$ for $1\leq k\ne i \leq  n-1$.
\elem

\bnota\label{main.n1}  Let the notations of  Lemma \ref{main.l3} hold.
Let $C^*$ be an open subset of $C$ containing a finite set $K$ and $C_i^*=\rho_{\alpha_i} (C^*)$.
Let $Q_0$ be as in Remark \ref{def.r1}, i.e.,   $Q_0\subset Q_{\rm reg}$, $\rho_{\alpha_i}|_{\rho_{\alpha_i}^{-1} (Q_0)} :\rho_{\alpha_i}^{-1} (Q_0)\to Q_0$ is a smooth morphism 
and $\delta_{\alpha_i}$ does not vanish on $\rho_{\alpha_i}^{-1} (Q_0)$.
Suppose that for every $i=1, \ldots, n-1$

(a) $C_i^*$ is contained in $Q_0$;

(b)  $\rho_{\alpha_i}|_{C^*}: C^* \to C_i^*$ is an isomorphism and the morphism $\rho_{\alpha_i}|_C : C \to Q$ is proper;

(c) for every $p \in K$ the vectors $\delta_{\alpha_1}(p), \ldots, \delta_{\alpha_{n-1}}(p)$ are linearly independent in $T_pX$ and
their span meets $T_pC$ at zero only.

Taking a smaller $C^*\supset K$ we can suppose that (c) holds for every $p \in  C^*$ as well
 (note that this implies
that  the vector fields $\delta_{\alpha_1}|_{C^*}, \ldots, \delta_{\alpha_{n-1}}|_{C^*}$
induce a basis of $N_XC^*$).
Choose functions $f_1, \ldots,  f_{i-1},f_{i+1}  \ldots  f_{n-1}\in \bk [Q]$ 
such that  they vanish on $C_i$ and at every point of $\rho_{\alpha_i} (K)$ their differential are linearly independent (this is possible since $C_i^*$ and $Q_0$ are smooth). 
Shrinking $C^*$ further we can suppose that the same is true for every point of $C_i^*$.
\enota

\blem\label{main.new.l1} Let Notation \ref{main.n1} hold. Then 
the matrix
 $T^i=\left[ \delta_{\alpha_k} (\rho_{\alpha_i}^*(f_j))\right]_{j,k=1; j,k \ne i}^{n-1}$ is invertible over $C^*$. 
\elem

\bproof  Let $p \in C^*$ and $q=\rho_{\alpha_i} (p) \in C_i^*$. Since $\rho_{\alpha_i}$ is smooth over $q$ by (a) in Notation \ref{main.n1},
$\delta_i(p)$ generates the kernel of $\rho_{\alpha_i*}: T_pX \to T_qQ$. Hence, the vectors 
$\rho_{\alpha_i*}(\delta_{\alpha_j}(p)), j\in \{ 1, \ldots, n-1 \}\setminus \{ i\}$ are linearly independent and their span meets $T_qC_i$ at zero only.
In particular, these vectors induce a basis in the fiber of the normal bundle $N_QC_i$ over $q$.
Since each $f_j$ vanishes on $C_i$, its differential $\dd f_j$ vanishes on $TC_i$. That is, the linear independence of
$\dd f_1, \ldots,  \dd f_{i-1}, \dd f_{i+1}  \ldots  \dd f_{n-1}$ implies that they induce a basis in the fiber of the conormal bundle of $C_i$ in $Q$
over $q$. This implies in turn that the value of the matrix $T^i$ at $p$ is invertible which yields the desired conclusion.
\eproof

\blem\label{main.l4} Let Notation \ref{main.n1} hold.
Suppose that $a$ is a nonzero regular function on $C$ that vanishes on $C\setminus C^*$ only, 
but with a sufficiently high multiplicity such that   for every $i$ the function $\frac{a}{ {\rm det} T^i}$ is regular on $C$.
Let $\bar b$ be any $(n-2)$-column vector whose entries are in $\bk[C]$.  Then  
there exists another $(n-2)$-column vector $\bar c_i$ with
entries in $\bk [C]$ such that $T^i \bar c_i =a \bar b$.  In particular, $T^i (a \bar c_i)= a^2 \bar b$.
\elem

\bproof Indeed, the matrix $a(T^i)^{-1} =a\frac{1}{\det T^i}{\rm adj} (T^i)$ (where ${\rm adj} (T^i)$ is the classical adjoint) is regular over $C$.
Hence, $\bar c_i =a(T^i)^{-1} \bar b$ yields the desired conclusion.
\eproof

To use Lemma \ref{main.l4} for automorphisms of $N_XC$  as in Lemma \ref{main.l3} (ii) we need the entries of $a\bar c_i$
to be pullbacks $\rho_{\alpha_i}^*(h_j)|_C$ of regular functions
$h_j|_{C_i}$ on $C_i$ which is guaranteed by the next fact.

\bprop\label{main.p1} { \rm (\cite[Theorem 2.3.1]{Ud})}  Let  $\tau : C \to \bar C$ be a proper birational morphism of affine curves such that
its restriction to an  open subset $C^*$ of $C$ is an isomorphism. Then there exists $m>0$ such that 
every regular function $a\in \bk[C]$ vanishing on $C\setminus C^*$ with multiplicity at least $m$ is a pullback of a regular function on $\bar C$.
\eprop

Recall that by Notation \ref{main.n1} $\rho_{\alpha_i}|_C : C \to Q$ is not only birational but also proper. Hence, Proposition \ref{main.p1} is applicable and
choosing the column $a^2\bar b$ in Lemma \ref{main.l4} so that  its $k_0$-entry is $a^2$ whereas the rest of them are zeros we find now
regular functions $h_1, \ldots ,h_{i-1},h_{i+1} \ldots,  h_{n-1} \in \bk[Q]$ for which one has $\sum_{j=1, j\ne i}^{n-1}(\rho_{\alpha_i}^*(h_j)\delta_{\alpha_{k_0}}(\rho_{\alpha_i}^*(f_j)))|_C=a^2$ and 
 $\sum_{j=1, j\ne i}^{n-1}(\rho_{\alpha_i}^*(h_j)\delta_{\alpha_{k}}(\rho_{\alpha_i}^*(f_j)))|_C=0$ for $k \ne k_0$.
Thus, we have the following.

\bprop\label{main.p2} Let Notation \ref{main.n1} hold. For every $1\leq k_0\ne i \leq n-1$ a locally nilpotent vector field $\sigma$ in Lemma \ref{main.l3} can be chosen so that
its flow induces an automorphism of $TX|_C$ such that $\breve \delta_{\alpha_k}\mapsto\breve \delta_{\alpha_k}$ for $1\leq k\ne k_0\leq n-1$ and
$\breve \delta_{\alpha_{k_0}}\mapsto \breve \delta_{\alpha_{k_0}}+t a^2d \breve \delta_{\alpha_i}$ where $t$ is a time parameter, $a$ is a fixed regular function on $C$ that 
vanishes on $C\setminus C^*$ only,  but with sufficiently high multiplicity and $d$ is any given regular function on $C$.  
In particular,  automorphisms of $TX|_C$ that have this form are $G$-induced.

\eprop

\bcor\label{main.c1} Let the notations of Proposition \ref{main.p2} hold.
Let  $\pr : TX|_C \to N_XC$ be the natural projection, $A=\bk[C]$ and $w_j=\pr (\breve \delta_{\alpha_j})$ for $j=1, \ldots , n-1$.
Suppose (without loss of generality) that $C\setminus C^*$ is the zero locus of a function $q\in \bk[C]$
(i.e., $\bk[C^*]=A[\frac{1}{q}]$).
Then $w_1, \ldots, w_{n-1}$ generate  a free $A$-module $M$ such that the module $M[\frac{1}{q}]=M\otimes_A A[\frac{1}{q}]$ coincides with $N_XC^*$. 
Furthermore, there is an automorphism $\alpha\in G$ of $X$ such that it induces an automorphism of $M$ given by
$w_k\mapsto w_k$ for $1\leq k\ne k_0\leq n-1$ and
$w_{k_0}\mapsto w_{k_0}+t a^2d w_i$. 
\ecor

\bproof The first statement was explained in Notations \ref{main.n1}.
The second one follows from an application of $\pr$ to the automorphism of $TX|_C$ in Proposition \ref{main.p2}.
\eproof

Letting the index $i$ in Proposition \ref{main.p2} run from $1$ to $n-1$ and taking a smaller neighborhood $C^*$ of $K$ (where $K$
is as in Notations \ref{main.n1}), if necessary, we get the next result.

\bcor\label{main.c2} Let  $M$ be as in Corollary \ref{main.c1} 
(i.e., every $A$-automorphism of $M$ can be presented as a matrix from ${\rm \Mat}_{n-1}(A))$.
Then $C^*\supset K$ can be chosen so that for some principal ideal $\fq=qA\subset A$ whose zero locus is  $C \setminus C^*$
every  element of $\FF_{n-1} (A, \fq)$ is $G$-induced (where  ${\rm \Mat}_{n-1}(A)$ and $\FF_{n-1} (A, \fq)$ are as in Notation \ref{bass.n1}).
\ecor

\proof[Proof of Theorem \ref{main.t1} (i)]  Recall that for general elements $\alpha_1, \ldots, \alpha_{n-1}$ of a perfect $G$-family $\cA$
the morphism $\rho_{\alpha_i}|_C : C \to Q$ is proper  for every $i$ 
because either $C$ is once-punctured or the conclusions of Proposition \ref{nor.p3} are valid by virtue of Convention \ref{nor.conv1}. 
Furthermore, by Propositions \ref{nor.p1} and \ref{nor.p3}    $\rho_{\alpha_i}|_{C^*}: C^* \to C_i^*$  is an isomorphism and $C_i^*$ is contained in $Q_0$
where $C^*\subset C$ is a neighborhood of a given finite subset $K$ of $C$. Thus, we have the assumptions (a) and (b) from Notations \ref{main.n1}.
The assumption (c) is also valid by virtue of Propositions \ref{nor.p1}(iii) and \ref{nor.p3}(4). Hence, we can use Corollaries \ref{main.c1} and \ref{main.c2}.

Suppose that $\bar v= (v_1, \ldots, v_{n-1})$ is a basis of $N_XC$ as a free $A$-module.
Let $\bar w= (w_1, \ldots , w_{n-1})$ where $w_i$ is as in Corollary \ref{main.c1}.
That is, $\bar w$ is a basis of a free $A$-submodule  $M$ of $N_XC$ as in Corollary \ref{main.c1}. 
Then for some matrix $D\in {\rm Mat}_{n-1} (A)$
we have $\bar w=D\bar v$. Let $q$ be as in Corollary \ref{main.c2}. Since the $A[\frac{1}{q}]$-module  $M[\frac{1}{q}]$ coincides with $N_XC^*=N_XC[\frac{1}{q}]$
by Corollary \ref{main.c1}, the matrix $D$ is invertible as an element of $ {\rm Mat}_{n-1} (A[\frac{1}{q}])$.
 Let $K'=C\setminus C^*$.
Then  by Proposition \ref{nor.p1} (resp. \ref{nor.p3}) and Corollary \ref{main.c1} we can find general elements 
$\alpha_1', \ldots,  \alpha_{n-1}'$ of some perfect $G$-family $\cA$ such that for the free $A$-module $M'$ with
the basis $\bar w' = (\pr(\breve \delta_{\alpha_1'}), \ldots, \pr (\breve \delta_{\alpha_{n-1}'}))$ there exists a  principal ideal $\fq'=q'A$ such that
its zero locus is disjoint from $K'$ and every  element of $\FF_{n-1} (A, \fq')$ is $G$-induced.
Furthermore, $\bar w'=D'\bar v$ for some matrix $D'$ similar to $D$. By construction, automorphisms of $N_XC$ whose matrices in
basis $\bar v$ belong to $\SL_{n-1}(A)\cap (D^{-1} \FF_{n-1}(A, \fq)D)$ or $\SL_{n-1}(A)\cap ((D')^{-1} \FF_{n-1}(A, \fq')D')$ are $G$-induced.
Hence,  by Proposition \ref{bass.p1} every element of $\EE_{n-1} (A)$ is $G$-induced.
Switching to dual spaces we see 
that every $A$-automorphism $\lambda :  \frac{I}{I^2} \to \frac{I}{I^2}$ that 
belongs to $E_{n-1} (A)$ is $G$-induced which concludes the proof. \hfill $\square$\\

\bnota\label{main.n2}  Let  $A$ be a commutative ring,  $I$
be the ideal in $A[u_1, \ldots, u_{n-1}]$ (where each $u_i$ is a variable) that consists of all polynomials whose constant terms are zero,
and $\bar u =( u_1, \ldots, u_{n-1})$.
Consider the space $\cP\subset I$ that consists of those polynomials whose degree is at most $k-1$
and the space $\tilde \cF$ of $(n-1)$-tuples  $\bar f = (f_1, \ldots , f_{n-1})$ of elements of $\cP$
(note that  modulo $I^k$ compositions of elements of $\tilde \cF$ are again elements of $\tilde \cF$).
 For every $f_i \in \cP$ one has $f_i =f_i^1+ \ldots + f_i^{k-1}$ where each $f_i^l$ is a homogeneous polynomial of degree $l$.
 Consider two $A$-linear maps $\chi : \tilde \cF \to \tilde \cF, \, \bar f \mapsto (f_1^1, \ldots, f_{n-1}^1)$ 
 and $\pi : \tilde \cF \to \tilde \cF, \, \bar f \mapsto (f_1^{k-1}, \ldots, f_{n-1}^{k-1})$, i.e., the image of $\pi$ coincides with the space $\cF$
 of $(n-1)$-tupes of $(k-1)$-homogeneous forms in $\bar u$ over $A$.
 Consider also the subset $\tilde \cF'$ of $\tilde \cF$  that consists of all  $\bar f \in \tilde \cF$ for which
$(\bar f -\pi (\bar f))(\bar u)=\bar u$. 
 \enota

\blem\label{main.l5}  Let $\bar f\in \tilde \cF'$ and $\bar g \in \tilde \cF$. Then

{\rm (1)} $\pi (\bar f \circ \bar g)= \pi (\bar g)+ \pi (\bar f) \circ \chi (\bar g)$ and
$\pi (\bar g \circ \bar f)=\chi (\bar g) \circ \pi (\bar f)+ \pi (\bar g)$
(in particular, if $\bar g\in \tilde \cF'$ one has $\pi (\bar f \circ \bar g)=\pi (\bar f) + \pi (\bar g)$);

{\rm (2)} $\bar f \circ \bar g - \pi ( \bar f \circ \bar g) = \bar g \circ \bar f - \pi ( \bar g \circ \bar f) =\bar g -\pi (\bar g)$.

Furthermore, let $\bar h \in \tilde \cF$ be such that $\bar h \circ \bar g (\bar u) =\bar u$ modulo $I^k$
(and, thus, $\chi (\bar h \circ \bar g) (\bar u) =\chi (\bar h )  \circ \chi (\bar g) (\bar u)=\bar u$).  Then

{\rm (3)} $\pi (\bar h \circ\bar f \circ \bar g)=\chi (\bar h) \circ \pi (\bar f) \circ \chi (\bar g)$ and $\bar h \circ\bar f \circ \bar  g\in \tilde \cF'$. 

\elem

\bproof Note that for every $\bar e \in \cF$ one has $\pi (\bar g \circ \bar e)=\chi (\bar g) \circ \bar e$ and
$\pi (\bar e \circ \bar g)= \bar e \circ \chi (\bar g)$. Since every element of $\tilde \cF'$ is of the form
$\bar f (\bar u)=\bar u +\bar e (\bar u)$ where $\bar e =\pi (\bar f)$, this implies (1) and (2).
Then by (2) and (1) we have $\bar f \circ \bar g= \bar g \circ \bar f - \pi (\bar g \circ \bar f) + \pi (\bar f \circ \bar g)
=\bar g \circ \bar f - \chi (\bar g) \circ \pi (\bar f) + \pi (\bar f )\circ \chi (\bar g)$. Applying (1) again to the composition of $\bar h$ with both sides of the
last formula we get $\pi (\bar h \circ\bar f \circ \bar g)=\pi (\bar f) - \chi (\bar h) \circ \chi (\bar g) \circ \pi (\bar f) +\chi (\bar h) \circ \pi (\bar f) \circ \chi (\bar g)
 =\chi (\bar h) \circ \pi (\bar f) \circ \chi (\bar g)$.   Note also that $\bar f \circ \bar g =\bar g \circ \bar f  +\bar e$ for some $\bar e \in \cF$.
 Hence, modulo $I^k$ one has $\bar h \circ\bar f \circ \bar  g =\bar f +\bar h \circ \bar e$. Since the last summand belongs to $\cF$, we
 get  $\bar h \circ\bar f \circ \bar  g\in \tilde \cF'$ which concludes the proof.
\eproof

\bnota\label{main.n3} 
Suppose that $C$ is a smooth curve in a normal quasi-affine variety $X$  such that $C$ is closed in an affine variety containing $X$.
Suppose also that $U\subset X_{\rm reg}$ is a neighborhood of $C$ in $X$ such that  
$C$ is a strict complete intersection in $U$ given by $u_1=\ldots =u_{n-1}=0$ where each $u_i $ belongs to the defining ideal $I$ of $C$ in $\bk [X]$.
That is, $\bar u=( u_1, \ldots, u_{n-1})$ can be viewed as a set of generators of the graded ring
$\frac{\bk [X]}{I^{k}}\simeq A\oplus \bigoplus_{l=2}^{k-1}\frac{I^{l-1}}{I^{l}}$ over $A=\bk [C]$.   
Then every  $A$-automorphism $\theta$ of $\frac{\bk[X]}{I^{k}}$  is given by its values on generators of
 $\frac{\bk[X]}{I^{k}}$ over $A$. 
That is, in the terminology of Notation \ref{main.n2} $\theta$ is given by a map of the form $\bar u \mapsto \bar f( \bar u)$ where $\bar f =(f_1, \ldots , f_{n-1}) \in \tilde \cF$ 
and the automorphism 
$\lambda : \frac{I}{I^2} \to \frac{I}{I^2}$ of the conormal bundle induced by $\theta$ is given by $\bar u \mapsto \chi (\bar f) (\bar u)$.
In particular, we have the natural map $\Phi :  \tilde H \to \tilde \cF$ from the group $\tilde H$ of $A$-automorphisms of  $\frac{\bk[X]}{I^k}$
which becomes a   semigroup homomorphism if one equips $\tilde \cF$ with the binary operation of composition modulo $I^k$.
   Let $\tilde \cF_0$ be the subset of $\tilde \cF$ consisting of all $\bar f$ for which the Jacobian $\left[ \frac{\p f_i}{\p u_j}\right] _{i,j=1}^{n-1}$ is equal to 1
modulo $I^{k-1}$. 
Note that $\chi (\bar f)$ can be viewed as an element of $\SL_{n-1} (A)$ for every $\bar f\in \tilde \cF_0$.
For $k\geq 3$ consider $\tilde \cF_0'=\tilde \cF_0\cap \tilde \cF'$. In particular, every $\bar f=(f_1, \ldots, f_{n-1})\in \tilde \cF_0'$ is of the form 
$\bar f(\bar u)=\bar u+\pi(\bar f) (\bar u)$.   It is a straightforward fact   observed in   \cite[proof of Lemma 4.13]{AFKKZ} that for such $\bar f$
 the assumption on the Jacobian is equivalent to the following  
 $$\div \pi(\bar f) =\frac{\p f_1^{k-1}}{\p u_1}+ \ldots + \frac{\p f_{n-1}^{k-1}}{\p u_{n-1}}=0.$$
 That is,  $\cF_0=\pi (\tilde \cF_0')$ is the space of $(n-1)$-tuples of $(k-1)$-forms  with zero divergence.
\enota

\blem\label{main.l6} Let Notations  \ref{main.n3} hold,  $G$ be a subgroup  of $\Aut (X)$ such
that $G$-induced $A$-automorphisms of  $\frac{\bk[X]}{I^l}$ have  Jacobian 1 (modulo $I^{l-1}$)  for every $l \geq 2$ and
let $H$ be the subgroup of  all $G$-induced $A$-automorphisms of  $\frac{\bk[X]}{I^k}$.
Suppose that $\Phi (H)\supset \tilde \cF_0'$ or, equivalently, $\cF^0:=\pi (\Phi (H)\cap \tilde \cF_0')$ coincides with $\cF_0$. Then  $\Phi (H)= \tilde \cF_0$.
\elem

\bproof We use induction on $k$ with the case of $k=2$ provided by Theorem \ref{main.t1}(i).  Furthermore, by construction an automorphism
$\alpha \in G$ whose restriction to $C$ induces $\lambda$ in Theorem \ref{main.t1}(i) is  a composition of flows
of locally nilpotent vector fields of the form $\rho_{\alpha_i}^*(hf) \delta_{\alpha_i}$ from Lemma \ref{main.l3}
where  $f,h\in \bk[Q]$ are such that $f$ vanishes on $\rho_{\alpha_i} (C)$. These vector fields satisfy the description of $\sigma$ in  Example \ref{main.exa1}.
Hence,  for every $l\geq 2$ the automorphism $\alpha$ induces    an $A$-automorphism $\frac{\bk[X]}{I^l}$  with Jacobian 1.    
Suppose that the statement is true for $k-1$. In particular, for every $\bar g \in \tilde \cF_0$ there exists an $\theta \in H$
such that $\bar g -\pi (\bar g)=\Phi (\theta) -\pi (\Phi (\theta))$. 
 Then $\bar f:=\Phi (\theta^{-1}) \circ \bar g={\rm id} +  \Phi (\theta^{-1}) \circ (\pi (\bar g) - \pi (\Phi (\theta)))$ is contained in $\tilde \cF'$
and has Jacobian 1 modulo $I^{k-1}$.  That is, $\bar f \in \tilde \cF_0'$
and, consequently,     $\Phi (\gamma ) =\bar f$ for some   $\gamma \in H$.   This implies that $\Phi (\theta \circ \gamma)=\bar g$
where $\theta \circ \gamma$ is induced by some $\alpha \in G$ . This yields the desired conclusion.     
\eproof

\blem\label{main.l6a}  Let the notations of Lemma \ref{main.l6} hold.

{\rm (i)} Then the map $\pi \circ \Phi : H \to  \cF_0$ is a group 
homomorphism (transforming multiplication to addition).

{\rm (ii)} There is a natural $\SL_{n-1} (A)$-action on $\cF$ such that $\cF^0$ is invariant with respect to this action.

\elem

\bproof Statement (i) follows from  Lemma \ref{main.l5} (1).  Let $\theta$ belong to $H$,
$\Phi (\theta)=\bar g$ and $ \Phi (\theta^{-1})=\bar h$. By Theorem \ref{main.t1} (i) $\chi \circ \Phi (H)$ is
isomorphic to $\SL_{n-1}(A)$, i.e., $\chi (\bar g)$ can be any element of $\SL_{n-1}(A)$ and $\chi (\bar h)$ is its inverse. 
Hence, one has the desired natural $\SL_{n-1}(A)$-action on $ \cF$
given by $\bar e \mapsto \chi (\bar h) \circ \bar e \circ \chi (\bar g)$. Let 
 $\pi (\bar f) \in \cF^0$ for some  $\bar f \in \tilde \cF_0'$, i.e., there exists
$\gamma \in H$ such that $\bar f =\Phi (\gamma) $. 
Since $\chi (\bar h) \circ \pi (\bar f) \circ \chi (\bar g) =\pi (\bar h \circ\bar f \circ \bar g)=\pi (\Phi (\theta^{-1} \gamma \theta))$
and $\bar h \circ\bar f \circ \bar  g\in \tilde \cF'$ by Lemma \ref{main.l5} (3),
we have $\chi (\bar h) \circ \pi (\bar f) \circ \chi (\bar g) \in  \cF^0$.  Therefore,
the $\SL_{n-1}(A)$-action sends every element of $\cF^0$ to an element of $\cF^0$  which yields (ii) and the desired conclusion. \eproof

We shall need also the following well-known fact from representation theory.

\bthm\label{main.t2} Let $F$ be the space of $(n-1)$-tuples of homogeneous $k$-forms in variables $\bar u=(u_1, \ldots , u_{n-1})$ over a field $\bK$
such that the divergence of every $q \in F$ is zero. Then $F$ is an irreducible $\SL_{n-1}(\bK)$-module under the
$\SL_{n-1} (\bK)$-action  given by $q(\bar u) \mapsto L^{-1} (q(L(\bar u))$ where $L \in \SL_{n-1} (\bK)$ and $q \in F$.
\ethm

\bproof  Denote by $V$ the $\bK$-vector space with the basis $\bar u$ and by $V^*$ its dual. 
Then the space of homogeneous $k$-forms in $\bar u$ can be treated
as ${\rm Sym}^k V$ and  the space $F$
as $V^* \otimes {\rm Sym}^k V$
(cf. \cite[Lemma 4.12(a) and Formula (13)]{AFKKZ}). 
In these notations divergence is the contraction map $V^* \otimes {\rm Sym}^k V \to {\rm Sym}^{k-1} V$ and $F$ is its kernel
which is invariant under the natural $\SL_{n-1}(\bk)$-action.  

Following \cite[page 216]{FuHa} we use the notations $\Gamma_{a_1, \ldots, a_{n-2}}$ (where $a_1, \ldots, a_{n-2}$ are nonnegative integers)
to denote irreducible representations of $\SL_{n-1}(\bk)$ classified by highest weights.  In particular,
$V^* \simeq \bigwedge^{n-2}V$ is the natural space of the irreducible representation $\Gamma_{0, \ldots ,0,1}$ of $\SL_{n-1}(\bk)$,
while for ${\rm Sym}^k V$
the corresponding irreducible representation is  $\Gamma_{k,0, \ldots ,0}$.
By \cite[Proposition 15.25 ]{FuHa}  ${\rm Sym}^k V \otimes \bigwedge^{n-2} V$ is the
sum of two irreducible representations  $\Gamma_{k-1, 0, \ldots , 0}$ and $\Gamma_{k, 0, \ldots , 0,1}$.
The first of them is, of course, the space ${\rm Sym}^{k-1} V$
and the second one is, therefore, the kernel of the contraction homomorphism. Hence, $F$ is irreducible and we are done. \eproof

\blem\label{main.l7} Let  Notation \ref{main.n3}  hold
(i.e., every automorphism of $\frac{I}{I^2}$ is given by a  matrix from $\GL_{n-1}(A)$) and let $G$
be generated by a saturated set $\cN \subset \LND (X)$.
Let all automorphisms of $\frac{I}{I^2}$ given by elements of  $\SL_{n-1}(A)$ be $G$-induced.
 Suppose that for every $x \in C$ there exists
a neighborhood $V \subset C$ and
$\sigma \in \cN$ with an associated partial quotient morphism $\tau : X\to P$
 such that $\sigma (x) \ne 0$, $\tau$ is smooth at $x$ and it
 induces an isomorphism $\tau|_V: V \to \tau (V)$ where
$\tau (V)$ is contained in the smooth part of $P$.
Then every $A$-automorphism $\theta : \frac{\bk[X]}{I^k} \to \frac{\bk[X]}{I^k}$   with Jacobian 1 is $G$-induced.
\elem

\bproof   By Lemma \ref{main.l6} it suffices to prove that $\cF_0=\cF^0$.
By Lemma \ref{main.l6a}  $\cF^0$ is a subgroup of $\cF_0$ which is invariant under the natural $\SL_{n-1}(A)$-action.
Hence, $\cF^0$ contains a unique maximal $A$-module $\cR$ which  is also invariant under the $\SL_{n-1}(A)$-action.
For every $A$-module $M$ denote by $M_{\mu}$ the localization of $M$ with respect to a maximal ideal $\mu\subset A$.
Recall that two submodules $M$ and $N$ of an $A$-module coincide if $M_{\mu}=N_{\mu}$ for every maximal ideal $\mu$ \cite[Proposition 3.9]{AM}.
Hence, by Lemma \ref{main.l6} it suffices to prove that $(\cF_0)_{\mu}=\cR_{\mu}$. 
By  \cite[Corollary 2.7]{AM} this is true if $\cR_{\mu} + \mu (\cF_0)_{\mu}=(\cF_0)_{\mu}$. 
In turn the latter equality is equivalent to the following $\frac{\cR_{\mu}}{\mu(\cF_0)_{\mu} \cap \cR_{\mu}}=\frac{(\cF_0)_{\mu}}{\mu(\cF_0)_{\mu}}$.
Note that $\frac{(\cF_0)_{\mu}}{\mu(\cF_0)_{\mu}}$ is the space of $(n-1)$-tuples of homogeneous $k$-forms in variables $\bar u=(u_1, \ldots , u_{n-1})$ over 
the residue field of $\mu$ such that the divergence is zero. 
By Theorem \ref{main.t2} it is enough to prove that  $\frac{\cR_{\mu}}{\mu(\cF_0)_{\mu} \cap \cR_{\mu}}$ is nonzero.
By the Nullstellensatz $\mu$ is the defining ideal of a point $x \in C$.  
Choose a regular function $f$ on $P$ that vanishes on $\tau (C)$ and has a nonzero differential at $\tau (x)$.
Then $\breve f=\tau^*(f)$ belongs to   the defining ideal $I$  of $C$ in $ \bk[X]$. 
Furthermore, since $\tau|_V: V \to \tau (V)$ is an isomorphism the assumptions on smoothness and the differential
imply that  $\breve f$ belongs to $I\setminus I^2$.
Consider the flow of the locally nilpotent vector field $\breve \sigma =\breve f^{k-1}\sigma$ with time parameter 1.
It sends $u_i$ to $u_i+\breve \sigma (u_i)+\sum_{l=2}^\infty \frac {\breve \sigma^l(u_i)}{l!}=u_i+\breve f^{k-1}\sigma (u_i) +\breve f^{2k-3}b$ where $b$
is some regular function on $X$. Note that if $k \geq 3$, then $\breve f^{2k-3}b$ belongs to $I^k$.
On the other hand, $\breve f^{k-1}\sigma (u_i) \in I^{k-1}$ and, furthermore, $\breve f^{k-1}\sigma (u_i) \in I^{k-1}\setminus \mu I^{k-1}$ provided that $\sigma (u_i)(x)\ne 0$.
Since $\sigma (x) \ne 0$ and $\sigma$ is not tangent to $C$ (because of the isomorphism $\tau|_V: V \to \tau (V)$) we do have $i$ for which  $\sigma (u_i)(x) \ne 0$.
Hence, $(\breve f^{k-1}u_1, \ldots, \breve f^{k-1}u_{n-1})$ yields  the desired nonzero element of 
$\frac{\cR_{\mu}}{\mu(\cF_0)_{\mu} \cap \cR_{\mu}}$ which concludes the proof.
\eproof

\proof[Proof of Theorem \ref{main.t1} (ii)]  By Propositions \ref{nor.p1} and \ref{nor.p3} letting $\sigma =\delta_{\alpha_1}$ we see that $\sigma$ satisfies the assumptions
of Lemma \ref{main.l7}. Since $\EE_{n-1}(A)=\SL_{n-1}(A)$ and by Theorem \ref{main.t1}(i) every element of $\EE_{n-1}(A)$ is 
$G$-induced, we see that all assumptions of Lemma \ref{main.l7} are valid.  This implies the desired conclusion. \hfill $\square$\\

\brem\label{main2.r1}  Let us emphasize that by  construction an automorphism
$\alpha \in G$ whose restriction to $C$ induces $\lambda$ (resp. $\theta$) in Theorem \ref{main.t1} is  a composition of  flows
of locally nilpotent vector fields of the form
$\rho_{\alpha_i}^*(hf) \delta_{\alpha_i}$ from Lemma \ref{main.l3}
where  $f,h\in \bk[Q]$ are such that $f$ vanishes on $C_i=\rho_{\alpha_i} (C)$.
Furthermore, this construction requires only the knowledge of  $f$ and $h$ on the $k$th infinitesimal neighborhood of $C_i$ in $Q$ and 
we do not care about the behavior of $f$ and $h$ on $Q\setminus C_i$.
Hence, one can require additionally that $f$ and $h$ vanish with a given multiplicity on a given closed subvariety of $Q$ disjoint from $C_i$. 
\erem

In conclusion of this section we want to pose the following problem.\\

{\bf Question}. {\em Does Theorem \ref{main.t1} remain true if the assumption $n\geq 4$ is replaced by $n\geq 3$?}

\section{Main Theorem. II}

The main aim of this section is to show that Theorem \ref{main.t1} is valid without the assumption that $C$ is irreducible.


\bthm\label{main2.t1} {\rm (1)} Theorem \ref{main.t1} is valid if assumption ($\bullet$) is replaced by one of the following:

{\rm (i)} either  every irreducible component of $C$ is once-punctured, or 

{\rm (ii)} Convention \ref{nor.conv1} holds.

{\rm (2)} Furthermore, let $l\geq 1$, $Z$ be a closed subvariety of $X$ disjoint from $C$, $\cN_{Z,l}$ be the subset of $\cN$ that consists of those elements of $\cN$
that vanish on $Z$ with multiplicity at least $l$ and  $G_{Z,l}$ be the subgroup of $G$ generated by the flows of elements of $\cN_{Z,l}$. 
Suppose that in (i) (resp. (ii))
${\rm codim}_X Z \geq 2$ (resp. ${\rm codim}_X Z \geq 3$).
Then  the automorphisms $\lambda$ and $\theta$  in Theorem \ref{main.t1} (modified as in (1)) are $G_{Z,l}$-induced. 
\ethm

\proof[Proof of Theorem \ref{main2.t1} (i).] Let  $C'$ be an irreducible component of 
$C$ and $C''=C\setminus C'$. Assume first that $C=C'$. The quasi-affine variety $X\setminus Z$ is $G_{Z,l}$-flexible
by \cite[Theorem 2.6]{FKZ}. Hence, the desired conclusion in this case follows from Theorem \ref{main.t1}.
If $C\ne C'$, then by induction on the number of irreducible components of $C$ one can suppose that the desired conclusion
holds for the curve $C''$.  That is, we have an automorphism $\alpha \in G_{Z,l}$ with a prescribed restriction to the $k$-infinitesimal 
neighborhood of $C''$. Replacing $Z$ by $Z \cup C''$ we get also an automorphism $\beta \in G_{Z\cup C'',l}\cap G_{C'',k}$
with a prescribed restriction to the $k$-infinitesimal neighborhood of $C'$. Hence, one can choose $\alpha$ and $\beta$
so that $\alpha \circ \beta \in G_{Z,l}$ has
a given restriction to the $k$-infinitesimal neighborhood of $C$ which concludes the proof. \hfill $\square$

\brem\label{main.2.r1} The above argument does not work for Theorem \ref{main2.t1} (ii), since  $X$ must be affine by
 Convention \ref{nor.conv1}. Hence, one cannot replace $X$ by $X\setminus Z$ (resp. $X\setminus C''$) when  
 $Z\ne \emptyset$ (resp. $C''\ne \emptyset$). 
\erem

Hence, we need the following technical result.

\blem\label{main2.l1} Let  $X$ be an affine $G$-flexible variety where $G$ is generated by a saturated set $\cN \subset \LND (X)$.
Let $\delta \in \cN$ and its partial quotient morphism $\rho : X \to Q$
be such that condition ($\#$) from Convention \ref{nor.conv1} is valid. Let $C$ be a closed curve in $X$
such that $C\subset X_{\rm reg}$ and $Z$ be a closed subvariety of $X$ with ${\rm codim}_X Z \geq 3$. 
Let $\cA$ be an algebraic $G$-family. 
Suppose that $W_\alpha=\rho_\alpha^{-1}(\overline {\rho_\alpha (Z)})$ where $\rho_\alpha =\rho \circ \alpha^{-1}$ for $\alpha \in \cA$
and  $\overline {\rho_\alpha (Z)}$ is the closure of ${\rho_\alpha (Z)}$ in $Q$. 
Then for every $l\geq 1$ there exists $\gamma \in G_{Z,l}$ such that
for a general $\alpha \in \cA$ the curve $\gamma (C)$ does not meet $W_\alpha$.
\elem

\bproof Let $\cX = \cA \times X$ and $\cZ=\cA \times Z$. Consider the morphism $\tilde \rho : \cX \to \cA \times Q, \, (\alpha, x) \mapsto (\alpha, \rho_\alpha (x))$.
 Let $\cV=\tilde \rho (\cZ)$, $\bar \cV$ be the closure of $\cV$ in $\cA\times Q$ and $\cD=\bar \cV \setminus \cV$.
By  Chevalley's theorem $\cV$ is a constructible set \cite[Chapter II, Exercise 3.19]{Har} and, hence, $\cD$ is a subvariety of $\bar \cV$.
Let $\kappa  : \bar \cV \to \cA$ be the natural projection, $\bar \cV'$ be an irreducible component of $\bar \cV$ 
and let $\cD'$ be an irreducible component of $\cD$ contained in $\bar \cV'$. 
By another Chevalley's theorem \cite[Chapter II, Exercise 3.22]{Har}
general fibers of  $\kappa|_{\bar \cV'}$ (resp. $\kappa|_{\cD'}$) are of the same  pure dimension
$\dim \bar \cV'- \dim \cA$ (resp. $\dim \bar \cD'- \dim \cA$). Since $\dim \cD' < \dim \bar \cV'$, we see that 
$\kappa^{-1}(\alpha) \cap \cD'$ cannot be an irreducible component of $\kappa^{-1}(\alpha) \cap \bar \cV'$ for a general $\alpha \in \cA$.
That is,  $\kappa^{-1}(\alpha) \cap \cD'$ is contained in the closure of $\kappa^{-1}(\alpha) \cap \cV$ in the fiber $\alpha \times Q$.
Hence, for a general $\alpha \in \cA$ the closure $\overline{\rho_\alpha (Z)}$ of $\rho_\alpha (Z)$
in $\alpha \times Q$ coincides with $\kappa^{-1} (\alpha) \cap \bar \cV$. In particular, for $\cW =\tilde \rho^{-1} (\bar \cV)=\tilde \rho^{-1} (\overline{ \tilde \rho (\cZ)}))$
and a general $\alpha \in \cA$ one has 
$\pi^{-1}(\alpha) \cap \cW= W_\alpha$ where $\pi : \cX \to \cA$ is the natural projection. 

 Let $G_{Z,l}$ act on $\cX$ by $\gamma. (\alpha, x) =(\alpha,  \gamma  (x))$.
Then this action preserves $\cZ$ and it is transitive on every fiber of $\pi|_{\cX \setminus \cZ}$ by \cite[Theorem 2.6]{FKZ}.
Let $\cC=\cA\times C$.
By Theorem \ref{pre.t1}  and \cite[Remark 1.9(1)]{Ka20} for a general $\gamma$ in some algebraic $G_{Z,1}$-family 
$\gamma.\cC$ meets $\cW$ along a subvariety $R$ of dimension at most
$\dim \cW + \dim \cC-\dim \cX=\dim \cA + \dim W_\alpha +1 -\dim X$ where $\alpha$ is a general element of $\cA$.
Note that by condition ($\#$)  $\dim W_\alpha \leq \dim X-2$ since $\dim Z \leq \dim X -3$.
Hence, $\dim R < \dim \cA$ and for a general $\alpha$ the set $R\cap \pi^{-1}(\alpha)=W_\alpha \cap \gamma (C)$ is empty which concludes the proof.
\eproof

\proof[Proof of Theorem \ref{main2.t1} (ii).]
Let $\cA$ be a family as in Proposition \ref{nor.p3} 
and $(\alpha_{1}, \ldots, \alpha_{n-1})$ be a general element of $\cA^{n-1}$.
That is, each $\rho_{\alpha_i}|_C : C\to  Q$ is a closed embedding and $C_i=\rho_{\alpha_i} (C)$ is contained in $Q_0$
where $Q_0$ is as in Remark \ref{def.r1}.  Replacing $C$ by $\gamma (C)$ where $\gamma \in G_{Z,l}$ is as in Lemma \ref{main2.l1}
we can suppose that $C$ does not meet $\rho_{\alpha_i}^{-1}( \overline{\rho_{\alpha_i} (Z)})$ for every $i=1, \ldots, n-1$.
This implies that $C_i$ does not meet the closed subvariety $\overline{\rho_{\alpha_i} (Z)}$ in $Q$.

Let  $C'$ be an irreducible component of 
$C$ and $C''=C\setminus C'$. Assume first that $C=C'$.  Recall that by  Remark \ref{main2.r1} 
an automorphism $\Psi \in G$ with a prescribed restriction $\theta$ to the $k$-infinitesimal neighborhood of $C$
can be obtained as a composition of automorphisms of the form
$\rho_{\alpha_i}^*(hf) \delta_{\alpha_i}$ where $h$ and $f$ are regular functions on $Q$ such that $f$ vanishes on $C_i$. 
Consider  $\tilde f$ and $\tilde h \in \bk [Q]$ such that  in the $k$th infinitesimal neighborhood 
of $C_i$ in $Q$ they coincide with $f$ and $h$ respectively but they   vanish on $\overline{\rho_{\alpha_i} (Z)}$ with multiplicity at least $l$
(we can find such functions since the closed subvarieties $C_i$ and $\overline{\rho_{\alpha_i} (Z)}$ of $Q$ are disjoint).
By Remark \ref{main2.r1} replacing each element $\rho_{\alpha_i}^*(hf) \delta_{\alpha_i}$  in the composition with $\rho_{\alpha_i}^*(\tilde h\tilde f) \delta_{\alpha_i}$ 
we get a new automorphism $\tilde \Psi\in G$ that induces the same automorphism $\theta$ of the $k$-infinitesimal neighborhood of $C$. 
However, $\tilde \Psi$ is already an element of $G_{Z, l}$ which concludes the case of irreducible $C$.

If $C\ne C'$, then by induction on the number of irreducible components of $C$ one can suppose that the desired conclusion
holds for the curve $C''$.  That is, we have an automorphism $\Phi'' \in G_{Z,l}$ with a prescribed restriction to the $k$-infinitesimal 
neighborhood of $C''$. Replacing $Z$ by $Z \cup C''$ we get also an automorphism $\Psi' \in G_{Z\cup C'',l}\cap G_{C'', k}$
with a prescribed restriction to the $k$-infinitesimal neighborhood of $C'$. Hence, $\Psi' \circ \Psi'' \in G_{Z,l}$ has
a prescribed restriction to the $k$-infinitesimal neighborhood of $C$ which concludes the proof. \hfill $\square$

\bcor\label{main2.c1} Let the assumptions of Theorem \ref{main2.t1} hold and Convention \ref{nor.conv1}  be valid. 
For every nonvanishing section $v$ of $N_XC$ there exists $\sigma \in \cN$ such that $v =\pr (\sigma)$ where $\pr : TX|_C \to N_XC$ is the
natural projection. In particular, $N_XC$ admits a basis
$ v_1 \ldots , v_{n-1}$ over $A$ such that each $v_i=\pr (\sigma_i)$  where $\sigma_1, \ldots , \sigma_{n-1}\in \cN$.
\ecor

\bproof  The second statement is especially simple. Indeed, by Proposition \ref{nor.p3}(5) we can suppose that $v_i=\pr (\sigma_i)$  where $\sigma_i \in \cN$ for $i\leq n-2$.
Then one can consider an $A$-automorphism $\theta$ of $N_XC$  given by $v_{n-1}\to v_1$, $v_1\to -v_{n-1}$ and $v_i\to v_i$ for $2\leq i \leq n-2$.
Since every element of $ \EE_{n-1}(A)$ is $G$-induced  by Theorem \ref{main2.t1} and $\theta \in \EE_{n-1} (A)$,
we see that $v_{n-1}=\pr (\sigma_{n-1})$ where $\sigma_{n-1}\in \cN$.

For the first statement note that $v=\sum_{i=1}^{n-1} a_iv_i$ where $a_1, \ldots , a_{n-1}\in A$ have no common zeros.
Since $n-1\geq 3$, by the Bass stable range theorem \cite[Theorem 11.1]{Ba} for some $b_i \in A, \, i\geq 1$ the ideal generated by 
$a_1', ,\ldots, a_{n-2}'$ is $A$ where $a_i'=a_i-b_ia_{n-1}, \, i=1, \ldots, n-2$.  Hence,  applying consequently elements of $\EE_{n-1}(A)$
we can send $v$ first to $\sum_{i=1}^{n-2}a_i'v_i$, then to $v_{n-1}+ \sum_{i=1}^{n-2}a_i'v_i$ and then to $v_{n-1}$. This yields the desired
conclusion.  \eproof

\section{Applications}

\bnota\label{stro.n1} 
In this section  $C_1$ and $C_2$ are smooth  curves
in a smooth quasi-affine variety $X$ with defining ideals $I_1$ and $I_2$ in the algebra $\bk [X]$ of regular functions on $X$. 
We suppose also that $C_1$ and $C_2$ are closed in an affine variety containing $X$.
Let $X$ possess 
a volume form $\omega$ (i.e., $\omega$ is a nonvanishing section of the canonical bundle on $X$) and let
each conormal bundle $\frac{I_j}{I_j^2}$ of $C_j$ in $X$ be trivial. 
Recall that by Lemma \ref{main.l2} there is a neighborhood $U_j$ of $C_j$ in $X$ in which  $C_j$ is a strict complete intersection 
given by $u_{1,j}=\ldots =u_{n-1,j}=0$ where  $u_{1,j}, \ldots , u_{n-1,j}\in I_j$. That is,
$\frac{\bk[X]}{I_j^k}\simeq A_j\oplus \bigoplus_{l=1}^{k-1}\frac{I_j^l}{I^{l+1}_j}$
is a graded algebra where  $A_j =\frac{\bk[X]}{I_j}$.
For $k \in \N$ we denote by
 $\theta : \frac{\bk[X]}{I_1^k} \to \frac{\bk[X]}{I_2^k}$ an isomorphism of these algebras. 
\enota

\bdefi\label{stro.d1} 
Since $N_X C_j$ is trivial, the existence of $\omega$ implies the existence  of a volume form on $C_j$.
Fix  volume forms $\omega_j$ on $C_j$ such that $\tilde \theta^* \omega_1=\omega_2$ where 
the isomorphism $\tilde \theta : C_2\to C_1$ is induced by 
$\theta$. Choose a section $\pr_j : TX|_{C_j}\to TC_j$ of the canonical inclusion $TC_j\hookrightarrow TX|_{C_j}$ 
and consider the section $\tilde \omega_j=\omega_j\circ \pr_j $ of the dual bundle $(TX|_{C_j})^\vee$ of $TX|_{C_j}$.
Then one can require that
 $\omega|_{C_j}$ coincides with $\tilde \omega_j \wedge \dd u_{1,j} \wedge \ldots \wedge \dd u_{n-1,j}$.
Under this requirement  the determinant of  the matrix $\left[\frac{\p \theta (u_{l,1})}{\p u_{m,2}}\right]_{l,m=1}^{n-1}$ 
is  well-defined modulo $I_2^{k-1}$ 
(i.e., it is independent of the choice of coordinates $u_{1,j}, \ldots , u_{n-1,j}$).
Hence, we say again that $\theta$ has Jacobian 1 if the determinant of  $\left[\frac{\p \theta (u_{l,1})}{\p u_{m,2}}\right]_{l,m=1}^{n-1}$ is equal to 1 modulo $I_2^{k-1}$.
\edefi

The next fact was is a straightforward  consequence of Theorem \ref{main2.t1}.

\bthm\label{stro.t1} Let the asumptions of Theorem \ref{main2.t1}  be satisfied with $X$ being a smooth quasi-affine variety equipped with a volume form.
Let $\Aut^1(X)\subset \Aut (X)$ be the subgroup of automorphisms of $X$ that have Jacobian 1 and let $Y$ be an affine variety containing $X$
as an open subset.
Suppose that  $\cF$ is a family of smooth curves in $X$ with trivial normal bundles such that each $C\in \cF$ is closed in $Y$ and
every  isomorphism between $ C_1\in \cF$ and $C_2\in \cF$ extends to an automorphism from $\Aut^1 (X)$ (resp. $\SAut (X)$). Then
every isomorphism $\theta : \frac{\bk[X]}{I_1^k} \to \frac{\bk[X]}{I_2^k}$ with the Jacobian 1 is $\Aut^1 (X)$-induced (resp. $\SAut (X)$-induced).
\ethm

\bproof    Let $\theta_1 : \frac{\bk[X]}{I_1} \to \frac{\bk[X]}{I_2}$ be the isomorphism generated by $\theta$. By the assumption $\theta_1^{-1}$ is generated by
some automorphism $\Phi \in G$ such that $\Phi (C_2)=C_1$. By Lemma \ref{main.l1} $\Phi$ generates an isomorphism $\varphi : \frac{A}{I_2^k} \to \frac{A}{I_1^k}$
with Jacobian 1. Hence, $\varphi \circ \theta : \frac{\bk[X]}{I_1^k} \to \frac{\bk[X]}{I_1^k}$ has Jacobian 1 and by Theorem \ref{main2.t1} it is generated 
by an automorphism $\Psi \in G$. Hence, $\theta$ is induced by $\Phi^{-1}\circ \Psi$ which is the desired conclusion.
\eproof

\brem\label{stro.r0} It is worth mentioning that in the case when $C_1=C_2$ we do not have to require  in Theorem \ref{stro.t1} that $\theta$ is an $A_1$-automorphism
as we did in Theorem \ref{main2.t1}.

\erem

\bexa\label{stro.exa1}  Let Notation \ref{stro.n1} hold and $\varphi : C_1 \to C_2$ be any isomorphism of smooth polynomial curves.

 (1) Consider the case of  $X=\SL_n (\bk), \, n\geq 3$. By van Santen's theorem  \cite{St} (see also \cite{Ka20}) there is an
automorphism $\alpha \in \Aut^1 (X)$ such that $\alpha|_{C_1}=\varphi$ and also $X$ has a volume form. 

(2)  More generally, let $X$ be  a connected complex algebraic group of dimension at least 4 without nontrivial characters.
The absence of characters implies that the natural action of $\SAut (X)$ on $X$ is transitive and, therefore, $X$ is a flexible
variety (\cite[Theorem 0.1]{AFKKZ}).  It has also a left invariant volume form. The Feller-van Santen theorem \cite{FS} states that
there is an automorphism $\alpha \in \Aut^1 (X)$ such that $\alpha|_{C_1}=\varphi$.\footnote{It is not stated explicitly in \cite{FS}
that the Jacobian of $\alpha$ is 1. However, by virtue of Lemma \ref{main.l1} this follows from the construction in \cite{FS}.}
 
  (3) Let $X$ be  an $(m-1)$-sphere over $\bk$ (i.e., a nonzero fiber of a non-degenerate quadratic form  $f$ on $\A^m$) where $m \geq 6$.
 Then $X$ is flexible and 
 by \cite[Theorem 8.3]{Ka20} there is an
automorphism $\alpha \in \SAut (X)$ such that $\alpha|_{C_1}=\varphi$.   By the adjunction formula there is a volume form $\omega$
on $X$ such that $\dd f\wedge \omega$ is the standard volume form on $\A^m$.

\eexa

Thus, we have the following.

\bcor\label{stro.c1} The conclusions of  Theorem \ref{stro.t1} are valid for all varieties in Example \ref{stro.exa1} (1)-(3)
in the case when $\cF$ is the family of smooth polynomial curves contained in $X$.

\ecor

\brem\label{stro.r1}
(i) The absence of nontrivial characters is essential for Example \ref{stro.exa1} (2). Without this assumption the statement is not true. 
Indeed, in the presence of nontrivial characters a linear algebraic group $X$ is isomorphic as an algebraic variety
to $X_0\times (\C^*)^m$ where $m>0$ and $X_0$ does not admit nonconstant invertible functions. 
In particular, the image $\pi (C_j)$ of a polynomial curve $C_j\subset X$ under the natural projection $\pi : X \to (\C^*)^m$ is a point $p_j$.
The absence of nonconstant invertible functions on $X_0$ implies that every automorphism
$\Phi$ of $X$ induces an automorphism $\phi$ of $(\C^*)^m$. In particular, $\pi_* \circ \Phi_*=\varphi_*\circ \pi_*$. However, one easily can find
an isomorphism $\lambda : N_XC_1\to N_XC_2$ for which $\pi_* \circ \lambda \ne \varphi_*\circ \pi_*$. In particular, such isomorphism
is not induced by an automorphism of $X$.

(ii) By virtue of the Lefschetz principle  Corollary \ref{stro.c1} as well as  the Feller-van Santen theorem are valid if one replaces the ground field $\C$ by
a universal domain $\bk$ \cite{Ek}.

\erem

Actually, for Example \ref{stro.exa1} (3) Corollary \ref{stro.c1} can be strengthened. To show this we need the next fact.

\bprop\label{stro.p2} Let $X$ be an $(m-1)$-sphere over $k$ where $m\geq 6$. 
Then $X$ is parametrically suitable. 
\eprop

\bproof  Parametric suitability follows immediately from Example \ref{prop.ex1} since $X$ is a homogeneous space of ${\rm SO}(m)$. 
\eproof

\brem\label{stro.r2} Example \ref{prop.ex1} and, therefore, the above proof are based on
a nontrivial result from \cite{DDK}. Hence, we give another rather trivial argument 
that implies weak parametric suitability which is sufficient for our purposes (see Remark \ref{prop.r1} (3)).
For the sake of notations we
consider the case $m=6$. Then one can
treat $X$ as a hypersurface in $\A^6$ whose coordinate form is $z_1z_2+z_3z_4+z_5z_6=1$.
 Consider the commuting locally nilpotent vector fields
$\delta_1 = z_4\frac{\p}{\p z_1}-z_2 \frac{\p}{\p z_3}$, $\delta_{2}=z_2 \frac{\p}{\p z_5} -z_6\frac{\p}{\p z_1}$
and  $\delta_{3}= z_6\frac{\p}{\p z_3}-z_4 \frac{\p}{\p z_5}$ on $X$. 
The kernel of $\sigma =a_1\delta_1+a_2\delta_2+a_3\delta_3, \, a_i \in \bk$
is generated by $z_2,z_4,z_6$ and $a_3z_1+a_2z_3+a_1z_5$ if not all $a_i$ are zero.
Choose the algebraic family $\cB\subset {\rm SO}(6)$
such that for every $\beta \in \cB$ the homomorphism $\beta^* : \bk^{[6]}\to \bk^{ [6]}$ preserves the subalgebra $\bk [z_2,z_4,z_6]$
and $(\beta^*)^{-1} (z_5)$ is of the form  $a_3z_1+a_2z_3+a_1z_5$. Then one can see that $\beta_* (\delta_1)$ has the same kernel
as $\sigma$. This implies that for $\delta=\delta_1$ and a general triple $(\beta, \gamma, \chi) \in \cB^3$ the span of 
$\Ker \beta_*(\delta) \cdot \Ker \gamma_*(\delta) \cdot \Ker \chi_*(\delta)$
contains all coordinates $z_1, \ldots , z_6$ and, thus, it
coincides with $\bk [X]$. Hence,  the morphism $(\rho_\beta, \rho_\gamma, \rho_\chi): X \to
Q_\beta\times Q_\gamma \times Q_\chi$ (where $\rho_\beta, \rho_\gamma$ and $\rho_\chi$ are the algebraic quotient morphisms) is a closed
embedding (and, in particular, proper) which is the desired conclusion.
\erem

\bthm\label{stro.t2} Let $X$ be
an $(m-1)$-sphere over $\bk$ where $m\geq 6$.  
Let $C_1$ and $C_2$ be isomorphic  smooth closed  curves in $X$ with  trivial normal bundles 
and $I_j$ be the defining ideal of $C_j$ in $\bk [X]$.
Suppose that $\EE_{m-2}(\bk[C_j]) =\SL_{m-2}(\bk[ C_j])$.
Then every automorphism $\theta : \frac{\bk[X]}{I^k} \to \frac{\bk[X]}{I^k}$ with Jacobian 1 is $\SAut (X)$-induced.
\ethm

\bproof 
As we mentioned in Example \ref{stro.exa1} (3)
$X$  is flexible, it has a volume form and every isomorphism of smooth closed curves in $X$ extends to an automorphism from $\SAut(X)$
by \cite[Corollary 8.4]{Ka20}.
The vector fields $\delta_1, \delta_2, \delta_3$ on $X$ from Remark \ref{stro.r2} have zero loci of
codimension 3 and algebraic quotient morphisms which are also smooth outside these zero loci. 
Hence, condition $(\#$) from Convention \ref{nor.conv1} is satisfied, while condition $(\#\#$) follows from Proposition \ref{stro.p2}.
Therefore, the assumptions of Theorem \ref{main2.t1} and, therefore, of Theorem \ref{stro.t1} are satisfied which yields the desired conclusion.
\eproof

\bcor\label{stro.c2} Let $X$ be an $(m-1)$-sphere over $\bk$ where $m\geq 6$.  
Let $C_1$ and $C_2$ be isomorphic  smooth closed  curves in $X$ such that their irreducible components are rational.
The every isomorphism of $k$th infinitesimal neighborhoods of $C_1$ and $C_2$ with Jacobian 1 extends
to an automorphism of $X$.
\ecor

\bproof By Theorem \ref{stro.t2} it suffices to show that $N_XC_j$ is a trivial vector bundle and $\EE_{m-2}(\bk[C_j]) =\SL_{m-2}(\bk[ C_j])$.
By  \cite{Se1} $N_XC_j$ is isomorphic to $T\oplus L$ where $T$ is a trivial vector bundle and $L$ is a line bundle on $C_j$.
Note that $L$ is trivial since the Picard group of a smooth affine rational curve is zero. Hence, $N_XC_j$ is trivial.
In the case when every irreducible component $C'$ of $C_j$ is a smooth polynomial curve the equality $\EE_{m-2}(\bk[C_j]) =\SL_{m-2}(\bk[ C_j])$ follows
from the fact that $\bk [C']$ is a Euclidean domain. Hence, in the general case of smooth rational irreducible components of $C_j$ the similar
equality follows from \cite[Chapter IV, Corollary 5.4]{Lam} and we are done.
\eproof

The next example demonstrates an application of Theorem \ref{main2.t1} (2).

\bexa\label{qua.exa1} Let $Z$ be a closed subvariety of $X=\A^n$ of dimension $\dim Z \leq n-3$ where $n\geq 4$.
Recall that $X\setminus Z$ is flexible by \cite{Wi}, \cite [Sec. 2.15,  p. 72, Exercise (b$'$)]{Gr} (see also \cite{FKZ})
and $X$ has a volume form inherited from $\A^n$. Therefore,
we have the assumptions and conclusions of Theorem \ref{main2.t1} (2) 
in this case. Furthermore,  if $n\geq 5$, then by \cite[Theorem 7.1]{Ka20} any isomorphism $\varphi : C_1\to C_2$ of two
smooth closed curves in $X$ disjoint from $Z$ extends to an automorphism of $X\setminus Z$ which is an element of $\SAut (X\setminus Z)$. 
Hence, arguing as in Theorem \ref{stro.t1} we see that if the normal bundles of $C_1$ and $C_2$ in $X$ are trivial,
then every isomorphism of the $k$th infinitesimal neighborhoods of $C_1$ and $C_2$ with Jacobian 1
extends to an automorphism from $\SAut (X\setminus Z)$.
\eexa

We want to mention one more example which was surprisingly unknown.

\bexa\label{stro.ex2}  Let $X \simeq \A^n=\A^2\times \A^{n-2}, \, n\geq 4$. 
Suppose that the first factor is equipped with a coordinate system $(z_1,z_2)$.
Let $\delta=\frac{\p}{\p z_1}$ and $\cB=\SL_2({\bk})$ act on the first factor, i.e.,  for $\beta\in \cB$ the vector field
$\delta_\beta=\beta_* (\delta)$ is of the form $a\frac{\p}{\p z_1}  +b\frac{\p}{\p z_2}$
and its algebraic quotient morphism $\rho_\beta : \A^n \to Q_\beta =\A^{n-1}$ is given by $(z_1, \ldots, z_n) \mapsto (bz_1-az_2, z_3, \ldots, z_n)$.
Then for general $\beta \ne \gamma \in \cB$ the morphism $(\rho_\beta, \rho_\gamma) : \A^n \to Q_\beta\times Q_\gamma$ is a closed embedding
and, thus, we have parametric suitability.
Every isomorphism of smooth closed curves extends to an automorphism of $X$  (\cite{Cr} and \cite{Je})
which belong to $\SAut (X)$. 
Hence, the assumptions and  the conclusions of Theorem \ref{stro.t1} are valid. It is also worth mentioning that the assumption that  the Jacobian is 1
can be replaced in this case by the assumption that the Jacobian is any nonzero constant.
\eexa

\section{Holomorphic case}

In this section we shall show that in the complex case the assumption that $\SL_{n-1}(A)=\EE_{n-1}(A)$ in Theorem \ref{main.t1} (ii)  (resp. \ref{main2.t1}) can be omitted
if we consider automorphisms of the normal bundle of $C$ induced by holomorphic automorphisms of $X$. 

\bdefi\label{hol.d1}
Let $\Hol (X)$ be the algebra of holomorphic functions on a complex affine algebraic variety $X$ and $C$ be a smooth closed curve in $X$
whose defining ideal  in $\Hol (X)$ is denoted by $\hI$. 
Suppose that for some holomorphic automorphism $\psi$ of $X$ one has $\psi (\hI) =\hI$. 
As before, $\psi$ generates an automorphism $\psi_k : \frac{\Hol (X)}{\hI^k} \to \frac{\Hol (X)}{\hI^k}$ which we call  holomorphically induced.
\edefi

The aim of this section is the following.

\bthm\label{hol.t1} Let 
$X$ be a normal complex affine $G$-flexible variety of dimension $n\geq 4$ where $G\subset \SAut (X)$
is a group generated by a saturated set $\cN\subset \LND (X)$,  $C$ be a smooth closed curve in $X$ such that $C\subset X_{\rm reg}$ and let $0\ne \delta\in \cN$ have 
a partial quotient morphism $\rho : X \to Q$ such that  Convention \ref{nor.conv1} is valid for this data.
 Let $\hI$ be the defining ideal of $C$ in  $\Hol (X)$ and $\hA =\frac{\Hol (X)}{\hI}$.
Suppose that  the normal bundle $N_XC$ of $C$ in $X$ is trivial.
Then  every $\hA$-automorphism $\theta : \frac{\Hol (X)}{\hI^k}\to \frac{\Hol (X)}{\hI^k}$   
 with Jacobian  1 is holomorphically induced.

\ethm

We start the proof with the following fact which is a simple consequence of the Ivarsson-Kutzschebauch theorem  \cite[Theorem 2.3]{IK}.

\bprop\label{hol.p1} Let $\hA$ be the ring of holomorphic functions on a smooth affine curve $C$. Then $\EE_n (\hA)=\SL_n(\hA)$
for every $n\geq 2$.
\eprop

\bproof
Consider
$\lambda \in \SL_{n} (\hA)$ as a holomorphic map $\theta : C \to \SL_{n}(\C )$.
Recall that such a map  is called null-homotopic if it can be continuously deformed to a constant map from $C$ into $\SL_{n} (\C)$.
As $C$ is a smooth affine curve, there exists a compact strong deformation retract $W\subset C$ such that $W=\bigvee_{i=1}^kW_i$
is a bouquet where each $W_i$ is homeomorphic to a circle.
Since $\SL_n(\C)$ is simply connected, for each $i$ there is a homotopy relative to the base points $H_i: W_i \times [0,1]\to \SL_n (\C)$
between $\theta|_{W_i} : W_i \to \SL_n(\C)$ and the constant map. Glueing the homotopies $H_i$ for $i=1, \ldots, k$
yields a homotopy relative to the base points of $\theta|_W: W \to \SL_n (\C)$
and the constant map. Hence, $\theta|_W$ and, therefore, $\theta$ is null-homotopic.
By the Ivarsson-Kutzschebauch theorem  \cite[Theorem 2.3]{IK} we have
$\SL_{n} (\hA) =E_{n} (\hA)$ which is the desired conclusion.
\eproof

We'll also need the following adjustment of Lemma \ref{main.l2}.

\blem\label{hol.l1} Let $Y$ and $Z$ be as in Lemma \ref{main.l2} and $Z^*$ be an open subset of $Z$
such that $Z\setminus Z^*$ is the zero locus of a function $g \in \bk [Z]$. Suppose that
$w_1, \ldots , w_{m-1}$ are sections of $N_YZ$ whose restrictions to $Z^*$ yield a basis of $N_YZ^*$ over $\bk [Z^*]$.
Then there exist functions $u_1, \ldots , u_{m-1}\in \bk [Y]$ that vanish on $Z$ and, thus, can be interpreted
as sections of the conormal bundle of $Z$ in $Y$, such that  $w_i(u_j)=g^k\delta_{ij}$ where $\delta_{ij}$ is the Kronecker symbol
and $k\geq 0$.
\elem

\bproof Extend $g$ to a regular function on $Y$ denoted by the same symbol. Then $Z^*$ is closed in the affine variety 
$Y^*=Y\setminus g^{-1}(0)$ and $w_1, \ldots , w_{m-1}$ form a basis of $N_{Y^*}Z^*$. Choose the dual basis 
in the conormal bundle which by Lemma \ref{main.l2} can be presented by regular functions $u_1', \ldots, u_{m-1}'\in \bk [Y^*]$.
For an appropriate $k \geq 0$ every function $u_i=u_i'g^k$ is regular on $Y$ which yields the desired conclusion.
\eproof

\bnota\label{hol.n1} Assuming Convention \ref{nor.conv1} let us fix notations for the rest of this section. Recall that by  Proposition \ref{nor.p3} there are general elements
$\alpha_1, \ldots, \alpha_{n-1}$ of a perfect $G$-family $\cA$ of automorphisms of $X$ such that every morphism $\rho_{\alpha_j}|_C : C \to Q$ is proper 
(where $\rho_{\alpha_j}: X \to Q$ is a partial quotient morphism of $\delta_{\alpha_j}$)  and it induces an isomorphism
$C\simeq C_j:=\rho_{\alpha_j}(C)\subset Q_0$. Letting $v_j =\pr (\delta_{\alpha_j}), \, j=1,\ldots , n-1$ we also note that by Corollary \ref{nor.c1}
there is a section $v_0$ of $N_XC$ such that $v_0,v_1, \ldots, v_{n-2}$ generate $N_XC$ as a free $A$-module where $A=\C [C]$.
In particular,  $v_{n-1} =q v_0+\sum_{i=1}^{n-2} a_i v_i$ where $q$ and each $a_i$ are in $A$.
Consider the $\hA$-module $M=N_XC\otimes_A\hA$ where  $\hA =\Hol (C)$. In particular, $M$ is a free $\hA$-module with basis
$\bar v=(v_0, v_1, \ldots, v_{n-2})$. Denote by $M'$ its free $\hA$-submodule with basis $\bar v' =(v_{n-1},  v_1, v_2,\ldots , v_{n-2})$.
Using these bases we identify $\hA$-automorphisms of $M$ and $M'$ with  matrices from $\GL_{n-1}(\hA)$. 
We also denote by ${\rm Ind} (M)$ (resp. ${\rm Ind}(M')$ those automorphisms of $M$ (resp. $M'$) that are induced
by holomorphic automorphisms of $X$.
\enota

\blem\label{hol.l2} Let $v_i^j=\kappa_j(v_i)$ where the map $\kappa_j : N_XC \to N_QC_j$
is induced by $\rho_{\alpha_j} : X \to Q$.
Then for every $j=1, \ldots, n-2$ the vector bundle $N_QC_j$ is generated as $\bk [C_j]$-module by
$v_0^j, v_1^j$, $\ldots, v_{j-1}^j, v_{j+1}^j, \ldots , v_{n-2}^j$. Furthermore, $ v_1^{n-1}, \ldots , v_{n-2}^{n-1}$
generate $N_QC_{n-1}^*$ where $C_{n-1}^*=\rho_{\alpha_{n-1}}(C\setminus S)$ and $S$ is the zero locus of $q$.
\elem

\bproof Let $x \in C$ and $x_j=\rho_{\alpha_j} (x)\in C_j$. For proving the first claim, it suffices to show that 
the vectors $v_0^j(x_j), v_1^j(x_j)$, $\ldots, v_{j-1}^j(x_j), v_{j+1}^j(x_j)$, $\ldots , v_{n-2}^j(x_j)$ are linearly independent.
Note that $\kappa_j \circ \pr = \pr_j\circ\rho_{\alpha_j*}$ where $\pr_j: TQ|_{C_j} \to N_QC_j$
is the natural projection.  Since the kernel of $\rho_{\alpha_j*}$ is generated by $\delta_{\alpha_j}$,
we see that the kernel of $\kappa_j$ is generated by $v_j$.
Hence,  dependence of the vectors above implies
dependence of  $ v_0(x), v_1(x), \ldots , v_j (x), \ldots , v_{n-2}(x)$ 
contrary to the fact
that $v_0, v_1, \ldots , v_{n-2}$ generate $N_XC$. Hence, we have the first claim.
Since  $v_{n-1} = qv_0+\sum_{i=1}^{n-2} a_i v_i$, we see that $v_1, \ldots , v_{n-1}$ form a basis of $N_X(C\setminus S)$.
Hence, the argument before implies that $ v_1^{n-1}, \ldots , v_{n-2}^{n-1}$ form a basis of $N_QC_{n-1}^*$ 
and we are done.
\eproof

\blem\label{hol.l3} Let $S=q^{-1}(0)\subset C$.  Suppose that the rows and columns of every matrix in $\GL_{n-1}(A)$ are enumerated
by indices $0,1, \ldots, n-2$ (instead of $1, \ldots, n-1$).
Then

{\rm (i)}  every matrix of the form $I+ce_{ji} \in \EE_{n-1} (\hA)$ where $0\leq i\leq n-2$ and  
$ j\in \{ 1, \ldots,  n-2\} \setminus \{ i\}$ is  contained in $ {\rm Ind} (M)$;

{\rm (ii)} $ {\rm Ind} (M)$ contains also all matrices of the form 
 $I+ a_lbe_{0k}  - a_kbe_{0l}\in \EE_{n-1}(\hA)$ where  $1\leq k\ne l \leq n-2$, $b \in \fq$ and $\fq\subset A$ is an ideal whose zero locus is $S$.
\elem

\bproof
By Lemmas \ref{hol.l1} and \ref{hol.l2} for every 
$1\leq j\leq n-2$ we can find
regular functions  $f_0, f_1$, $\ldots, f_{j-1}, f_{j+1}, \ldots , f_{n-2}\in \C[Q]$ vanishing on $C_j$
such that  $v_i^j( f_k) =\delta_{ik}$ where $\delta_{ik}$ is the Kronecker symbol and  $k\ne i \in \{ 0, \ldots, n-2\} \setminus \{ j \} $.  
Recall that by Corollary \ref{main2.c1} for every $i=0, \ldots, n-1$ there exists a locally nilpotent vector field $\sigma_i \in \LND (X)$
such that $v_i=\pr (\sigma_i)$ (e.g., for $1\leq i \leq n-1$ one can put $\sigma_i=\delta_{\alpha_i}$).
We also have the equality $\kappa_j \circ \pr = \pr_j\circ\rho_{\alpha_j*}$ where $\pr_j: TQ|_{C_j} \to N_QC_j$ and $\kappa_j : N_XC \to N_QC_j$
are the natural projections.  Thus, $v_i^j=\pr_j\circ \rho_{\alpha_j*}(\sigma_i)$.
Since the differential $\dd f_k$ vanishes on $TC_j$, one has 
\be\label{hol.eq1} \delta_{ik} = v_i^j(f_k)=\rho_{\alpha_j*}(\sigma_i)(f_k)|_{C_j} =\sigma_i(\rho_{\alpha_j}^*(f_k))|_{C}.\ee
Suppose that $h$ is any holomorphic function on $Q$. 
Consider the complete vector field $\rho_{\alpha_j}^*(hf_i) \sigma_j$.
By Formula \eqref{pre.eq3} in Lemma \ref{pre.l1} and Formula \eqref{hol.eq1} its flow generates $\hA$-automorphisms of $M$
such that $v_k\mapsto v_k$ for $k \in \{ 0, \ldots, n-2\} \setminus \{ i \} $ and 
$v_i \mapsto v_i +t \rho_{\alpha_j}^*(h) |_C \, v_j, \, t \in \C$. 
This yields (i).

Similarly, Lemmas \ref{hol.l1} and \ref{hol.l2} imply that
one can find some $m> 0$ and
regular functions $f_1, \ldots , f_{n-2}$ on $Q$
such that $v_{i}^{n-1}(f_k) =q^m\delta_{ik}$ for $1\leq i \ne k \leq  n-2$.
Formulas \eqref{pre.eq3} and \eqref{hol.eq1} imply again that
every automorphism of $M'$ which is of the form $I+cq^me_{0j}, \, c \in \hA$ is contained in ${\rm Ind }(M')$.
Note that one has the following relation $\bar v'=D \bar v$ between the bases of $M$ and $M'$ where up to the change of indices $D$ is 
of the same form as the transpose of the matrix $\sigma$ in Lemma \ref{bass.l2}.
Since $n\geq 4$ and $\SL_{n-2}(\hA)=\EE_{n-2}(\hA)$  by Proposition \ref{hol.p1}, we have (ii)  by Proposition \ref{bass.p2} which concludes the proof.
\eproof

\blem\label{hol.l4} Let the assumptions of Lemma \ref{hol.l3} hold. Then $\EE_{n-1}(\hA)$ is
contained in ${\rm Ind}(M)$.
\elem

\bproof
Recall that $v_{i}=\pr (\delta_{\alpha_{i}})$ where $\alpha_1, \ldots, \alpha_{n-1}$ are general elements of some perfect $G$-family $\cA$.
 While keeping $\alpha_1, \ldots, \alpha_{n-2}$ intact, choose another general element $\alpha_{n-1}'\in \cA$ and let $v_{n-1}'=\pr (\delta_{\alpha_{n-1}'})$. 
 By Corollary \ref{nor.c1} we have $v_{n-1}' = q'v_0+\sum_{i=1}^{n-2} a_i' v_i$.  By Lemma \ref{hol.l3} this yields  holomorphically induced automorphisms
 of the form 
 $\theta'=I+ a_l'b'e_{0k}  - a_k'b'e_{0l}$ where  $1\leq k\ne l \leq n-2$ and $b'\in \hA$ is divisible by $(q')^m$ for some large $m>0$.
  By Proposition \ref{nor.p2}(4) we can suppose that 
 at some point $x\in C$ the vector $\delta_{\alpha_{n-1}'} (x)$ is such that the function $a_la_k'-a_l'a_k$ does not vanish at $x$
 and, therefore, $a_la_k'-a_l'a_k$ is a nonzero function. Let $b=ca_k'(qq')^m$ and $b'=-ca_k(qq')^m$ where $c \in \hA$. Then the holomorphically induced
 automorphism $\theta \theta'$ is given by $I+c(a_la_k'-a_l'a_k)(qq')^m e_{0k}$.
Thus, ${\rm Ind} (M)$ contains all automorphisms given by $I+de_{0k}$ where $d$ belongs to a nonzero ideal $\fp\subset \hA$.
Choose now new general elements $\alpha_{n-1}$ and $\alpha_{n-1}'$ in $\cA$. By Proposition \ref{nor.p2}(4)
we can suppose that for every $x$ in the zero locus $K$ of $\fp$ the vectors $\delta_{\alpha_{n-1}} (x)$ and $\delta_{\alpha_{n-1}'} (x)$
are such that the new function $(a_la_k'-a_l'a_k)(qq')^m$ does not vanish at $x$. Hence,  ${\rm Ind} (M)$
contains every automorphism  of the form $I+de_{0k}, \, d \in \fq$ where $\fq\subset \hA$ is an ideal whose zero locus in $C$ is disjoint from
$K$. By the weak Nullstellensatz \cite{Car}  $\fp +\fq =\hA$ and, hence, every elementary transformation
of the form $I+ae_{0k}, \, a \in A$ is in ${\rm Ind} (M)$. In combination with  Lemma \ref{hol.l3} this implies  the desired conclusion.
\eproof

\proof[Proof of Theorem \ref{hol.t1}] 
By Proposition \ref{hol.p1}  and Lemma \ref{hol.l4} we see now that ${\rm Ind} (M)$ contains $\SL_n(\hA)$.  It remains to note that  Lemma \ref{main.l7} 
remains valid if one replaces $A$ by $\hA$ and $G$ by the group of holomorphic automorphisms of $X$. 
Hence, this  holomorphic analogue of   Lemma \ref{main.l7} implies that
every  $\hA$-automorphism $\theta : \frac{\Hol (X)}{\hI^k}\to \frac{\Hol (X)}{\hI^k}$   
with Jacobian  1 is holomorphically induced and we are done.
\hfill $\square$\\

\brem\label{hol.r1} The fact that $\hA$ is not a Dedekind domain is the reason why in the proof of Theorem \ref{hol.t1} instead of Theorems \ref{bass.t1} and \ref{bass.t2}
we use  Proposition \ref{bass.p2}.
\erem

Repeating the argument in the proof of Corollary \ref{main2.c1} we get the following.

\bcor\label{hol.c1}  Let the assumptions of Theorem \ref{hol.t1} hold. Suppose that
$M=N_XC\otimes_A\hA$ is as in Notation \ref{hol.n1} (i.e., $M$ is a holomorphic normal bundle of $C$ in $X$).
Then every nonvanishing section of $M$ is induced by a complete holomorphic vector field.
\ecor

\bcor\label{hol.c2} Let $X$ be isomorphic to a connected  complex linear algebraic group  without nontrivial characters
and $\dim X  \geq 4$. Suppose that $C$ is a smooth closed curve in $X$ with a trivial normal bundle
and  defining ideal $\hI$ in the algebra $\Hol (X)$ of holomorphic functions on $X$. 
Let $\hA =\frac{\Hol (X)}{\hI}$ be the algebra of holomorphic functions on $C$.
Then  every $\hA$-automorphism $\frac{\Hol (X)}{\hI^k}\to \frac{\Hol (X)}{\hI^k}$   
 with Jacobian  1 extends to a holomorphic automorphism of $X$.

\ecor

\bproof  
By Theorem \ref{hol.t1} it suffices to check  the validity of Convention \ref{nor.conv1}, i.e., conditions ($\#$) and ($\# \#$) of that convention.
Condition  ($\#$) is obvious since multiplication produces free $\G_a$-actions on $X$. 
Condition ($\# \#$) follows from the followng fact: every connected  linear algebraic group $Y$ of dimension $n \geq 2$ without nontrivial
characters is parametrically suitable. Indeed, by Mostow's theorem \cite{Mo} $Y$ is isomorphic (as a variety) to the product $S\times U$ where $U$ 
is unipotent and $S$ is reductive. The absence of nontrivial characters implies that $S$ is semi-simple. If $S$ is not trivial, then we
have parametric suitability by Theorem \ref{prop.t1}.
Otherwise, $Y \simeq \A^n$ where $n \geq 2$ and the argument in Example \ref{stro.ex2} works which concludes the proof.
\eproof

\bcor\label{hol.c3} Let $X$ be isomorphic to $\SL_n (\C)$ where $n\geq 3$. Suppose that $C_1$ and $C_2$ are  isomorphic smooth closed curves in $X$
with defining ideals $I_1$ and $I_2$ in $\bk [X]$.  Suppose also that the normal bundles $N_XC_j, \, j=1,2$ are trivial.
Then every isomorphism $\theta : \frac{\bk[X]}{I_1^k} \to \frac{\bk[X]}{I_2^k}$ with  Jacobian 1 is holomorphically induced.
\ecor

\bproof  
By \cite[Corollary 11.12]{Ka20} every isomorphism $\varphi : C_1\to C_2$ extends to a holomorphic isomorphism of $X$
which by construction has Jacobian 1. Hence, it suffices to consider the case when $\varphi$ is the identity map
which is a special case of Corollary \ref{hol.c2}. Hence, we are done.
\eproof

\end{document}